\documentclass{gtart}

\def\ifplaintex{\expandafter\ifx\csname documentclass\endcsname\relax}

\ifplaintex 
\else
\fi

\def\gtm{{\mathsurround=0pt\it $\cal G\mskip-2mu$eometry \&\ 
$\cal T\!\!$opology $\cal M\mskip-1mu$onographs}}    

\def\gtp{{\mathsurround=0pt\it $\cal G\mskip-2mu$eometry \&\ 
$\cal T\!\!$opology $\cal P\!$ublications}}  

\def\recd{{\small Received:\qua\receiveddate\ifx\reviseddate\relax
\else\qquad Revised:\qua\reviseddate\fi\par}} 

\def\volumenumber#1{\def\thevolumenumber{#1}}
\def\volumeyear#1{\def\thevolumeyear{#1}}
\def\volumename#1{\def\thevolumename{#1}}
\def\papernumber#1{\def\thepapernumber{#1}}
\def\pagenumbers#1#2{\def\startpage{#1}\def\finishpage{#2}}
\def\published#1{\def\publishdate{#1}}
\def\received#1{\def\receiveddate{#1}}
\def\revised#1{\def\reviseddate{#1}}
\def\accepted#1{\def\accepteddate{#1}}

\long\def\asciiabstract#1{\long\def\theasciiabstract{#1}}

\let\\\par
\let\thevolumenumber\relax\let\thepapernumber\relax
\let\thevolumeyear\relax\let\startpage\relax
\let\finishpage\relax\let\publishdate\relax\let\receiveddate\relax
\let\reviseddate\relax\let\accepteddate\relax\let\theasciititle\relax
\let\theasciiauthors\relax
\let\theasciiabstract\relax

\let\theerratum\relax\let\theasciiemail\relax
\let\theshortauthors\relax\let\theshorttitle\relax

\def\startpage{1}\def\finishpage{15}\def\thepapernumber{77}

\volumenumber{2}
\volumename{Proceedings of the Kirbyfest}
\volumeyear{1999}

\long\def\maketitlep{   

\count0=\startpage

\gtm\nl        
{\small Volume \thevolumenumber: \thevolumename\nl 
\ifx\theerratum\relax\else Erratum \erratumnumber\nl\fi
Pages \startpage--\finishpage\nl}

\vglue 0.1truein   

{\parskip=0pt\leftskip 0pt plus 1fil\def\\{\par\smallskip}{\ifplaintex\large
\else\Large\fi\bf\thetitle}\par\medskip}   
\vglue 0.05truein 

{\parskip=0pt\leftskip 0pt plus 1fil\def\\{\par}{\sc\theauthors}
\par\medskip}%
 
\vglue 0.03truein 

{\small\leftskip 25pt\rightskip 25pt{\bf Abstract}\stdspace\theabstract

{\bf AMS Classification}\stdspace\theprimaryclass
\ifx\thesecondaryclass\relax\else; \thesecondaryclass\fi\par
{\bf Keywords}\stdspace \thekeywords\par}\vglue 7pt

}   

\font\phead=cmsl9 scaled 950
\font=cmsl9 scaled 1050
\font\pnum=cmbx10 scaled 913
\font\lnum=cmbx10 
\font\pfoot=cmsl9 scaled 950
\font=cmsl9 scaled 1050
\ifplaintex
\headline{\vbox to 0pt{\vskip -4.5mm\line{\small\phead\ifnum
\count0=\startpage ISSN 1464-8997 (on line)
1464-8989 (printed) \hfill {\pnum\folio}\else\ifodd\count0\def\\{ }%
\ifx\theshorttitle\relax\thetitle\else\theshorttitle\fi\hfill{\pnum\folio}
\else\def\\{ and }{\pnum\folio}\hfill\ifx\theshortauthors\relax\theauthors
\else\theshortauthors\fi\fi\fi}\vss}}
\footline{\vbox to 0pt{\vglue 0mm\line{\small\pfoot\ifnum\count0=\startpage
Published \publishdate:\qua\copyright\ \gtp\hfill\else
\gtm, Volume \thevolumenumber\ (\thevolumeyear)\hfill\fi}\vss
}}
\else
\makeatletter
\def\@oddhead{{\small\lhead\ifnum\count0=\startpage ISSN 1464-8997 (on line)
1464-8989 (printed) \hfill {\lnum\number\count0}\else\ifodd\count0
\def\\{ }\ifx\theshorttitle\relax \thetitle \else\theshorttitle\fi\hfill
{\lnum\number\count0}\else\def\\{ and }{\lnum\number\count0}
\hfill\ifx\theshortauthors\relax 
\theauthors\else\theshortauthors\fi\fi\fi}}\def\@evenhead{@oddhead}
\def\@oddfoot{\small\lfoot\ifnum\count0=\startpage Published \publishdate:\qua\copyright\ \gtp\hfill\else
\gtm, Volume \thevolumenumber\ (\thevolumeyear)\hfill\fi}
\def\@evenfoot{@oddfoot}
\makeatother
\fi

\let\maketitlepage\maketitlep
\let\makeshorttitle\maketitlepage
\let\maketitle\maketitlepage

\newwrite\gtoutfile
\long\gdef\makeheadfile{  
{\def\\{, }\def\s{ }
\immediate\openout\gtoutfile head.xxx
\immediate\write\gtoutfile{To: math@arxiv.org}
\immediate\write\gtoutfile{Subject: put OR rep NNNNN:ppppp}
\immediate\write\gtoutfile{--text follows this line--}
\immediate\write\gtoutfile{Proxy-for: \ifx\theasciiauthors\relax
\theauthors\else\theasciiauthors\fi\s<\ifx\theasciiemail\relax\theemail\else\theasciiemail\fi>}
\immediate\write\gtoutfile{\noexpand\\}
\immediate\write\gtoutfile{Authors: \ifx\theasciiauthors\relax
\theauthors\else\theasciiauthors\fi}
{\def\\{ }\immediate\write\gtoutfile{Title: \ifx\theasciititle\relax
\thetitle\else\theasciititle\fi}}
\immediate\write\gtoutfile{Subj-class: GT or SG, GR etc}
\immediate\write\gtoutfile{MSC-class: \theprimaryclass\ifx\thesecondaryclass\relax\else, \thesecondaryclass\fi}
\immediate\write\gtoutfile{Journal-ref: Geom. Topol. Monogr. \thevolumenumber\s
(\thevolumeyear) \startpage-\finishpage}
\immediate\write\gtoutfile{Comments: Published by Geometry and Topology Monographs at}
\immediate\write\gtoutfile{\s\s\s  http://www.maths.warwick.ac.uk/gt/GTMon\thevolumenumber/paper\thepapernumber.abs.html}
\immediate\write\gtoutfile{\noexpand\\}
\immediate\write\gtoutfile{}
\ifx\theasciiabstract\relax
\immediate\write\gtoutfile{\theabstract}\else
\immediate\write\gtoutfile{\theasciiabstract}\fi
\immediate\write\gtoutfile{}
\immediate\write\gtoutfile{\noexpand\\}
\immediate\write\gtoutfile{}
\immediate\closeout\gtoutfile}}  

\def\maketitlepage{\maketitlep\makeheadfile}
\let\makeshorttitle\maketitlepage
\let\maketitle\maketitlepage

\volumenumber{4}
\volumename{Invariants of knots and 3-manifolds (Kyoto 2001)}
\volumeyear{2002}
\papernumber{11}
\pagenumbers{161}{181}
\received{19 December 2001}
\revised{6 August 2002}
\accepted{10 September 2002}
\published{21 September 2002}

\usepackage{amsmath,amssymb,epsf}

\def\ostar{\star\hspace{-0.35cm} \circlearrowleft}
\def\ZLMO{\mbox{Z}^{LMO}}
\def\la{\langle\hspace{-0.05cm}\langle}
\def\ra{\rangle\hspace{-0.05cm}\rangle}
\def\ZHS{{\mathbb Z}HS^3}
\def\JS{{\cal A}(\star_{\{k_1,\ldots,k_\mu\}})}
\def\TS{{\cal A}(\star_{\{k\}})}
\def\JSS{\mbox{${\cal A}(\ostar_{\{k_1,\ldots,k_\mu\}})$}}
\def\Qset{{\mathbb Q}}
\def\Zset{{\mathbb Z}}
\def\WT{{\mbox{Tr}^{\circlearrowleft}}}
\def\thetdiag{
\begin{picture}(1.5,1)(-0.2,-0.3)
\qbezier[10](0,0)(0,0.5)(0.5,0.5)
\qbezier[10](1,0)(1,0.5)(0.5,0.5)
\qbezier[10](0,0)(0,-0.5)(0.5,-0.5)
\qbezier[10](1,0)(1,-0.5)(0.5,-0.5)
\qbezier[10](0,0)(0.5,0)(1,0)
\end{picture}
}
\def\thetdiagm{
\begin{picture}(1.5,0.5)(0.2,-0.15)
\put(0.70,0.62){\small $p(k)$}
\put(0.70,-0.27){\small $q(k)$}
\put(0.85,0.5){\vector(1,0){0.3}} 
\put(0.85,0){\vector(1,0){0.3}}
\qbezier[10](0.5,0)(0.5,0.5)(1,0.5)
\qbezier[10](1.5,0)(1.5,0.5)(1,0.5)
\qbezier[10](0.5,0)(0.5,-0.65)(1,-0.65)
\qbezier[10](1.5,0)(1.5,-0.65)(1,-0.65)
\qbezier[10](0.5,0)(1,0)(1.5,0)
\end{picture}
}
\def\thetdiagmz{
\begin{picture}(1.5,0.5)(0.2,-0.15)
\put(0.8,0.5){\vector(1,0){0.4}}
\put(0.8,0){\vector(1,0){0.4}}
\put(0.8,-0.5){\vector(1,0){0.4}}
\qbezier[10](0.5,0)(0.5,0.5)(1,0.5)
\qbezier[10](1.5,0)(1.5,0.5)(1,0.5)
\qbezier[10](0.5,0)(0.5,-0.5)(1,-0.5)
\qbezier[10](1.5,0)(1.5,-0.5)(1,-0.5)
\qbezier[10](0.5,0)(1,0)(1.5,0)
\end{picture}
}
\def\thetdiagsm{
\begin{picture}(0.75,0.5)(-0.2,-0.1)
\qbezier[10](0,0)(0,0.25)(0.25,0.25)
\qbezier[10](0.5,0)(0.5,0.25)(0.25,0.25)
\qbezier[10](0,0)(0,-0.25)(0.25,-0.25)
\qbezier[10](0.5,0)(0.5,-0.25)(0.25,-0.25)
\qbezier[10](0,0)(0.25,0)(0.5,0)
\end{picture}
}

\newtheorem{thm}{Theorem}
\newtheorem{lem}[thm]{Lemma}
\newtheorem{cor}[thm]{Corollary}
\theoremstyle{definition}
\newtheorem{defn}[thm]{Definition}
\newtheorem{notn}[thm]{Notation}

\newtheorem*{rems}{Remarks}

\newenvironment{fake}{\relax}{\relax} 

\begin{document}

\title
[A surgery formula for the 2-loop piece]
{A surgery formula for the 2-loop piece\\of the LMO invariant of a pair}
\author{Andrew Kricker}
\address{Department of Mathematics, University of Toronto\\Ontario, 
M5S 1A1, Canada}
\email{akricker@math.toronto.edu}

\begin{abstract}
Let $\Theta(M,K)$ denote 
the 2-loop piece of (the logarithm of) the LMO invariant of a knot $K$ in $M$,
a ${\mathbb Z}HS^3$. 
Forgetting the knot (by which we mean setting diagrams with legs to zero)
specialises $\Theta(M,K)$ to $\lambda(M)$, 
Casson's invariant.
This note describes an extension of Casson's surgery formula for his 
invariant to $\Theta(M,K)$. To be precise, we describe the effect on
$\Theta(M,K)$ of a surgery on a knot which together with $K$ forms a
boundary link in $M$. Whilst the presented 
formula does not characterise $\Theta(M,K)$,
it does allow some insight into the underlying topology.
\end{abstract}

\asciiabstract{Let \Theta(M,K) denote the 2-loop piece of (the
logarithm of) the LMO invariant of a knot K in M, a ZHS^3.  Forgetting
the knot (by which we mean setting diagrams with legs to zero)
specialises \Theta(M,K) to \lambda(M), Casson's invariant.  This note
describes an extension of Casson's surgery formula for his invariant
to \Theta(M,K). To be precise, we describe the effect on \Theta(M,K)
of a surgery on a knot which together with K forms a boundary link in
M. Whilst the presented formula does not characterise \Theta(M,K), it
does allow some insight into the underlying topology.}


\begin{fake} \end{fake}
\makeshorttitle

\setcounter{section}{-1}

\section{Introduction}

The simplest characterisation of $\lambda$,
Casson's invariant of integral homology three-spheres, is the following.
Let $K$ be an $f$-framed knot, where $f$ is plus or minus 1, in $M$, a
$\ZHS$; and let $M_K$ denote the result of surgery on $K$. Then,
\begin{equation}\label{casson}
\lambda(M_K) = \lambda(M) + f\, a_2(M,K).
\end{equation}
In the equation above 
$a_2(M,K)$ denotes
the coefficient of $k^2$ in the power series
$A_{(M,K)}(e^{k})$, where 
$A_{(M,K)}(t)$ is the symmetric Alexander polynomial of a knot $K$
in $M$, a $\ZHS$, (by symmetric is meant the representative satisfying
$A_{(M,K)}(t^{-1}) = A_{(M,K)}(t)$
and $A_{(M,K)}(1) = 1$).

What actually 
happened, of course, is that this formula 
was used to ``discover'' Casson's invariant 
in the image of the LMO invariant.
But let us reverse this history and ask how a hypothetical
student of finite-type invariants, ignorant of Casson's formula,
might ``derive'' it from the LMO invariant
(we recall the identification in detail in Section 2). 

To begin, (considering $S^3_K$ for simplicity), 
Casson's invariant is the 
co-efficient of the theta graph in the image of the LMO invariant:
\begin{equation}\label{LMO}
\ZLMO(S^3_K) = 1 + \frac{\lambda(S^3_K)}{2}\thetdiag
+ \mbox{terms of higher degree}.
\end{equation}
The Alexander polynomial, on the other hand, is to be found in the image
of the Kontsevich invariant, amongst the coefficients of the ``wheel'' 
diagrams (see Theorem \ref{abelianthm}):
\setlength{\unitlength}{1cm}
\begin{equation}\label{asubtwo}
\hat{Z}(K) = \hat{Z}(U)\left( 1 -\frac{1}{2} a_2(S^3,K) 
\ 
\begin{picture}(1.2,0.5)(0.9,0.15)
\qbezier[6](1.25,0.75)(1.25,1)(1.5,1)
\qbezier[6](1.25,0.75)(1.25,0.5)(1.5,0.5)
\qbezier[6](1.75,0.75)(1.75,1)(1.5,1)
\qbezier[6](1.75,0.75)(1.75,0.5)(1.5,0.5)
\qbezier[6](1.25,0.75)(1,0.75)(1,0.5)
\qbezier[6](1,0.5)(1,0.375)(1.125,0.25)
\qbezier[6](1.125,0.25)(1.25,0.125)(1.25,0)
\qbezier[6](1.75,0.75)(2,0.75)(2,0.5)
\qbezier[6](2,0.5)(2,0.375)(1.875,0.25)
\qbezier[6](1.875,0.25)(1.75,0.125)(1.75,0)
\end{picture}\ + \mbox{terms of higher degree}
\right).
\end{equation}

According to the LMO surgery formula, 
$Z^{LMO}(S^3_K)$ is obtained from the Kontsevich integral of $K$
by ``gluing chords'' into its legs, in a certain way
(see Equation \ref{LMOSURG}).
Casson's formula then follows from the obvservation that the
only term in the surgery formula which leads to a theta diagram 
is when a single chord is glued 
into the legs of the wheel with 2 legs:
\begin{equation}\label{pair}
\left<
-\frac{1}{2} f
\begin{picture}(1,0.7)(-0.25,-0.2)
\qbezier[16](0,-0.1)(0,0.25)(0.25,0.25)
\qbezier[16](0.25,0.25)(0.5,0.25)(0.5,-0.1)
\end{picture} \ ,-\frac{1}{2}a_2(S^3,K)\
\begin{picture}(1.2,0.5)(0.9,0.15)
\qbezier[6](1.25,0.75)(1.25,1)(1.5,1)
\qbezier[6](1.25,0.75)(1.25,0.5)(1.5,0.5)
\qbezier[6](1.75,0.75)(1.75,1)(1.5,1)
\qbezier[6](1.75,0.75)(1.75,0.5)(1.5,0.5)
\qbezier[6](1.25,0.75)(1,0.75)(1,0.5)
\qbezier[6](1,0.5)(1,0.375)(1.125,0.25)
\qbezier[6](1.125,0.25)(1.25,0.125)(1.25,0)
\qbezier[6](1.75,0.75)(2,0.75)(2,0.5)
\qbezier[6](2,0.5)(2,0.375)(1.875,0.25)
\qbezier[6](1.875,0.25)(1.75,0.125)(1.75,0)
\end{picture}\ \right> = f a_2(S^3,K)\  \frac{1}{2}\hspace{-0.15cm}
\thetdiagsm\ \ .
\end{equation}
The simple aim of this article is to tell this familiar story, and then to
repeat it, but in a more general context. This more general
story concerns $\Theta(M,K)$,
the 2-loop piece of the LMO invariant of a pair. We may think of
it as a Casson's invariant of $M$ equipped with a knot $K\subset M$.
Indeed, if we set all diagrams 
with legs to zero we recover $\frac{\lambda(M)}{2}\theta$. 

More formally: 
in this setting
diagrams may have their edges labelled by power series in a single variable,
\hspace{0.15cm} \vspace{-0.15cm}\begin{math}
\begin{picture}(1,0.5)(0,0)
\qbezier[20](0,-0.25)(0,0)(0,0.25)
\qbezier[6](0,0.25)(0,0.5)(-0.1,0.5)
\qbezier[6](0,-0.25)(0,-0.5)(-0.1,-0.5)
\put(0,0){\circle*{0.1}}
\put(0,-0.15){\vector(0,1){0.4}}
\put(0.2,0){$a_0 + a_1 k_1 + a_2 k^2 + \ldots\ \ . $ }
\end{picture}
\end{math}\hspace{3.75cm}\vspace{0.5cm}
This describes a series of Jacobi
diagrams: add a leg for every
power of $k$
(see Notation \ref{legs}).  In these terms, $\Theta(M,K)$ is defined
to be that part of the LMO invariant arising from marked thetas
(supressing labels below):
\begin{eqnarray*}
\ZLMO(M,K) & = & (\mbox{wheels})\,\sqcup\, \\
& & \ \ \ \ \ 
\mbox{exp}_{\sqcup}(\ \ \sum\thetdiagmz \ + \mbox{conn. diags. with $\geq$ 3
loops}\ ).
\end{eqnarray*}

At first glance, $\Theta(M,K)$ may not seem an interesting
invariant. It appears to be a rather arbitrary collection of diagrams
from the image of the LMO invariant. But here is the thing:

\begin{thm}[\cite{K}, Rozansky's conjecture]
$\Theta(M,K)$ is rational. That is, it is expressible as 
a finite combination of labelled theta diagrams, where each
label is a quotient of the form $\frac{p(e^k)}{A_{(M,K)}(e^k)}$,
for $p(e^k)$ a polynomial in $e^k$.
\end{thm}


Our main theorem (Theorem \ref{mainth}) 
concerns the effect on $\Theta(M,K)$ of 
the following move on $K$: $\pm 1$-framed
surgery on a knot $K'$ which together with
$K$ forms a boundary link. (See \cite{GK} for 
a theory of the rationality of the Kontsevich integral of
a knot or a boundary link.)
It observes a
generalisation of Casson's formula (Equation \ref{casson})
of the following general form. The supressed labels are of the
form $p(e^k)/A_{(M,K)}(e^k)$, $p(e^k)$ a polynomial, and the sum
is finite. 
\begin{equation}
\Theta( (M,(K,K'))_{K'} ) = 
\Theta( M , K ) + \sum \thetdiagmz.
\end{equation}

The contributing terms arise from the LMO surgery formula
in a familiar way (see Section 3): 
\setlength{\unitlength}{1.1cm}
\begin{equation}
\left<
-\frac{f}{2}
\begin{picture}(1,0.7)(-0.25,-0.2)
\qbezier[16](0,-0.1)(0,0.25)(0.25,0.25)
\qbezier[16](0.25,0.25)(0.5,0.25)(0.5,-0.1)
\end{picture} \ ,\
\begin{picture}(1.2,0.5)(0.9,0.15)
\put(1.375,1){\vector(1,0){0.3}}
\put(1.2,1.12){\small $p(k)$}
\put(1.2,0.22){\small $q(k)$}
\put(1.375,0.5){\vector(1,0){0.3}}
\qbezier[6](1.25,0.75)(1.25,1)(1.5,1)
\qbezier[6](1.25,0.75)(1.25,0.5)(1.5,0.5)
\qbezier[6](1.75,0.75)(1.75,1)(1.5,1)
\qbezier[6](1.75,0.75)(1.75,0.5)(1.5,0.5)
\qbezier[6](1.25,0.75)(1,0.75)(1,0.5)
\qbezier[6](1,0.5)(1,0.25)(1.125,0)
\qbezier[6](1.75,0.75)(2,0.75)(2,0.5)
\qbezier[6](2,0.5)(2,0.25)(1.875,0)
\end{picture}\ \right> = -f \thetdiagm\ \ .
\end{equation}

Section 1, ``The 1-loop piece of $\ZLMO(M,L)$ for a boundary link $L$'',
introduces the background for the second term above, i.e., the
wheel with 2 legs with marked edges. This arises from a certain 
non-commutative generalisation of the Alexander polynomial of a knot to
a 2-component boundary link.

Section 2, ``Casson's invariant'', 
recalls the identification of the Casson invariant 
with the 2-loop piece of the LMO invariant of a $\Zset HS^3$.
The proof of our generalised formula will be an adaptation of this
proof. Then, Section 3 recalls some elements of the recent theory of
the {\it rational expansion of the Kontsevich invariant}. This will describe
the origin and character of $\Theta(M,K)$. 

Our extension of Equation \ref{casson} to 
$\Theta(M,K)$ (see Theorem \ref{mainth}) is the subject of Section 4,
``Surgery on a sublink of a boundary link''.
Whilst this is our advertised goal,
our ulterior motive is to use this discussion 
to highlight, in the {\bf simplest possible terms},
one or two results from the recent 
theory of the rationality of the Kontsevich invariant, 
and expose some techniques from these (sometimes dense) papers.

We illustrate this technique in Section 5, where a step is taken
towards finding a Seifert surface based formula for $\Theta(M,K)$.

\medskip
{\bf Acknowledgements}\qua The author was partially 
supported by a Golda Meir fellowship,
and would like to thank Dror Bar-Natan for his support.
He would also like to thank
Tomotada Ohtsuki and Hitoshi Murakami for the truly memorable
workshop ``Invariants of knots and 3-manifolds''; and
also Stavros Garoufalidis and Tomotada Ohtsuki for discussions.

\section{The 1-loop piece of $Z^{LMO}(M,L)$ for a boundary\break link
$L$}

\paragraph{The diagram-valued determinant}
Consider Equation \ref{asubtwo}: actually, 
one can interpret the entire 
1-loop piece of the Kontsevich invariant of a knot (that is, 
the projection to the subspace generated by the {\it wheels}),
as a certain representation of the Alexander polynomial of a knot. 
To do this, write down
the usual definition, but replace the determinant you have written 
with the {\it diagram-valued}
determinant, to follow.

Consider, then, the (usual) determinant $|M|$ 
of a square matrix $M$ of power series in a single variable $k$
which augments to a matrix $M_{\varepsilon}$
which is invertible over $\Qset$ (one {\it augments}
$M$ by setting $k$ to zero). Such a determinant satisfies the 
following equation (which, observe, is well-defined, $(1-MM_\varepsilon^{-1})$
being small in the $k$-adic topology):
\begin{equation}\label{usualdet}
|M| = |M_{\varepsilon}|\, \mbox{exp}\left( -\mbox{Tr}\left(
\sum_{l=1}^{\infty} \frac{(1-MM_{\varepsilon}^{-1})^l}{l} \right) \right).  
\end{equation}
This equality follows (for example)
from a few straightforward
{\bf determinant-like properties} of the right-hand side $\Psi'(M)$, 
which we take this opportunity to collect, below. 
(In detail: DLP's 1 and 2 
imply that $\Psi'(M)$ is
unaffected by elementary row operations. Having assumed that $M$ 
augments to an invertible matrix, elementary row operations 
may be used to transform
$M$ to upper-triangular form, whence the result follows from 
DLP's 3 and 4.)  
All matrices mentioned below 
augment to invertible matrices. 
\begin{enumerate}
\item{$\Psi'(M) = |M|$ if $M=M_\varepsilon$.}
\item{$\Psi'(M_1 M_2) = \Psi'(M_1) \Psi'(M_2)$.}
\item{If $M$ is of the form
$
\left[
\begin{array}{ll}
A & B \\ 0 & C
\end{array}
\right],
$
then $\Psi'(M) = \Psi'(A)\Psi'(C)$.}
\item{$\Psi'([x]) = x, x\in \Qset[[k]].$}
\end{enumerate}
For example, property (2) follows from the BCH formula together with
the cyclic invariance of the trace (which kills commutators). 

The diagram-valued determinant, to be introduced presently,
is naturally motivated by the problem: Define a ``determinant'' 
$\Psi(M)$ of a square matrix $M$ 
of elements of $\Qset\la k_1,\ldots,k_\mu\ra $ (where $M$ 
augments to an invertible matrix). The notation
$\Qset\la k_1,\ldots,k_\mu \ra$ denotes the ring of power
series in $\mu$ non-commuting variables.
To be ``precise'': 
Define a ``useful''
function of such matrices of non-commutative power series 
satisfying DLP's one through four.

We might try to use the right-hand side of Equation \ref{usualdet} to
{\it define} such a determinant. As it stands, however, the expression
is inadequate: $Tr(AB)$ is no longer equal to $Tr(BA)$ if the
entries of $A$ and $B$ do not commute, and as a consequence
multiplicativity (property (2)) fails.

We can, however, restore
the cyclic invariance of the trace
by replacing the trace operation
with a certain {\it wheel-valued} trace. To recall this,
a notation for certain series of Jacobi diagrams helps.
By Jacobi diagram we mean the familiar thing (unitrivalent diagrams modulo 
AS, STU, IHX relations). This note, for example, will frequently use
a space 
${\cal A}(\star_{\{k,k'\}})$, 
the diagrams of which have univalent vertices (no skeleton)
labelled by either $k$ or $k'$. 
The object we require we will call
a {\it generating diagram}: a trivalent graph whose edges may have some
extra oriented (i.e. the pair of incident edges ordered)
bivalent vertices each of which has an element of $\Qset\la k_1,\ldots,k_\mu
\ra $ 
affixed to it. Such a diagram denotes an element of $\JS$, which we
will say {\it corresponds to} (or is {\it generated by})
the generating diagram. For example:

\begin{notn}\label{legs}
\setlength{\unitlength}{0.9cm}
\begin{eqnarray*} 
&
\begin{picture}(1,1)(0,0)
\qbezier[20](0,-0.5)(0,0)(0,0.5)
\qbezier[6](0,0.5)(0,0.75)(-0.1,0.75)
\qbezier[6](0,-0.5)(0,-0.75)(-0.1,-0.75)
\put(0,0){\circle*{0.1}}
\put(0,-0.15){\vector(0,1){0.4}}
\put(0.2,0){$(a_0 + a_1 k_1 + a_2 k_2 + a_3 k_1 k_2 + a_4 k_2 k_1 + \ldots 
\,$)}
\end{picture} \hspace{5cm} &  \\
& & \nonumber \\
= & a_0 
\begin{picture}(1,1)(-0.25,0)
\qbezier[20](0,-0.5)(0,0)(0,0.5)
\qbezier[6](0,0.5)(0,0.75)(-0.1,0.75)
\qbezier[6](0,-0.5)(0,-0.75)(-0.1,-0.75)
\end{picture} 
\hspace{-0.5cm} 
+\ \ 
a_1
\begin{picture}(1,1)(-0.25,0)
\qbezier[20](0,-0.5)(0,0)(0,0.5)
\qbezier[6](0,0.5)(0,0.75)(-0.1,0.75)
\qbezier[6](0,-0.5)(0,-0.75)(-0.1,-0.75)
\qbezier[12](0,0)(0.25,0)(0.5,0)
\put(0.65,-0.1){$k_1$}
\end{picture} 
\ \ \ \ +\ \ 
a_2
\begin{picture}(1,1)(-0.25,0)
\qbezier[20](0,-0.5)(0,0)(0,0.5)
\qbezier[6](0,0.5)(0,0.75)(-0.1,0.75)
\qbezier[6](0,-0.5)(0,-0.75)(-0.1,-0.75)
\qbezier[12](0,0)(0.25,0)(0.5,0)
\put(0.65,-0.1){$k_2$}
\end{picture} 
\ \ \ \ +\ \ 
a_3
\begin{picture}(1,1)(-0.25,0)
\qbezier[20](0,-0.5)(0,0)(0,0.5)
\qbezier[6](0,0.5)(0,0.75)(-0.1,0.75)
\qbezier[6](0,-0.5)(0,-0.75)(-0.1,-0.75)
\qbezier[12](0,0.25)(0.25,0.25)(0.5,0.25)
\qbezier[12](0,-0.25)(0.25,-0.25)(0.5,-0.25)
\put(0.65,0.15){$k_2$}
\put(0.65,-0.35){$k_1$}
\end{picture} 
\ \ \ \ +\ \ 
a_4  
\begin{picture}(1,1)(-0.25,0)
\qbezier[20](0,-0.5)(0,0)(0,0.5)
\qbezier[6](0,0.5)(0,0.75)(-0.1,0.75)
\qbezier[6](0,-0.5)(0,-0.75)(-0.1,-0.75)
\qbezier[12](0,0.25)(0.25,0.25)(0.5,0.25)
\qbezier[12](0,-0.25)(0.25,-0.25)(0.5,-0.25)
\put(0.65,0.15){$k_1$}
\put(0.65,-0.35){$k_2$}
\end{picture} 
\ \ \ \ +\ \ \ldots
& \nonumber \\
& & \nonumber
\end{eqnarray*}
\end{notn}

\begin{defn}

\

\begin{enumerate}
\item{
The {\it wheel-valued trace},
$\WT$, takes an element of ${\cal M}(\Qset\la k_1,\ldots, k_\mu \ra)$,
the set of square matrices
of elements of $\Qset\la k_1,\ldots,k_\mu \ra$,
to
\setlength{\unitlength}{0.5cm}
\[
\sum_{i}\ 
\begin{picture}(1,1)(-0.5,-0.2)
\qbezier[12](-0.5,0)(-0.5,1)(0,1)
\qbezier[12](0.5,0)(0.5,1)(0,1)
\qbezier[12](-0.5,0)(-0.5,-1)(0,-1)
\qbezier[12](0.5,0)(0.5,-1)(0,-1)
\put(0.5,0){\circle*{0.2}}
\put(0.5,-0.2){\vector(0,1){0.6}}
\put(0.75,-0.2){$M_{ii}$}
\end{picture}\hspace{0.75cm} \in \JS.
\]
}
\item{
Let
${\cal M}(\Qset\la k_1,\ldots,k_\mu \ra)^{\varepsilon}
\subset
{\cal M}(\Qset\la k_1,\ldots,k_\mu \ra)$,
denote the subset of matrices which
augment to invertible matrices. The {\it diagram-valued determinant}
$\Psi\tilde{\,}(M): {\cal M}(\Qset\la k_1,\ldots,k_\mu \ra)^{\varepsilon}
\rightarrow 
\JS
$ 
is defined by
\begin{equation}
M \mapsto |M_\varepsilon|\,
\mbox{exp}_{\sqcup}\left(
- \WT \left(
\sum_{l=1}^{\infty} \frac{(1-MM_{\varepsilon}^{-1})^l}{l} \right) \right)
\ 
\end{equation}
}
\end{enumerate}
\end{defn}
\begin{lem}
$\Psi\tilde{\,}$ satsifies determinant-like properties (1) through (4).
\end{lem}
The discussion below is in terms of the normalisation $\Psi(M) = 
\frac{1}{|M_\varepsilon|}\Psi\tilde{\,}(M)$.
\paragraph{The 1-loop piece}
(Note on terminology: by ``1-loop piece'' we mean the 1-loop piece
of the logarithm of $\ZLMO(M,K)$.)

The story which leads to the following theorem
begins with a certain influential conjecture and proof -- the conjecture
due to Melvin and Morton \cite{MM}, and the proof, one of the first
successes of the theory of finite-type invariants, due
to Bar-Natan and Garoufalidis \cite{BG}. 

Our theorem below connects the Kontsevich integral of a boundary link
to a certain abelian invariant of boundary links. By abelian we mean
arising from classical techniques, depending ultimately on certain linking
numbers. By boundary link we mean a link for the components of
which we can find a
a set of disjoint Seifert surfaces.

So, let $L$ be a boundary link of $\mu$ 0-framed
components in $M$, a $\Zset HS^3$. 
Take a set of 
disjoint Seifert surfaces $\Sigma=\{\Sigma_1,\ldots,\Sigma_\mu\}$ 
for the components of $L$. 
Say $\Sigma_j$ has genus $g_j$. On each  
surface $\Sigma_j$ choose a system $\{c^j_1,\ldots,c^j_{2g_j}\}$
of oriented curves which present a basis for
the first 
homology of $\Sigma_j$. Let $S$ denote the Seifert matrix corresponding
to these choices; that is, $S$ is a square matrix with rows and columns
indexed by the ordered list $\{c^1_1,\ldots,c^1_{2g_1},\ldots,
c^\mu_{1},\ldots,c^\mu_{2g_\mu} \}$. Observe that the Seifert matrices
of the individual surfaces appear in blocks along the diagonal, and that
the off-diagonal blocks measures 
linking between curves from different surfaces. We require, furthermore,
the following matrix of power series of $\Qset\la k_1,\ldots,k_\mu \ra$,
which we sometimes write $T(k_1,\ldots,k_{\mu})$:
\[
T = \mbox{Diag}( \underbrace{e^{k_1},\ldots,e^{k_1}}_{2g_1},
\ldots,
\underbrace{e^{k_\mu},\ldots,e^{k_\mu}}_{2g_\mu}).
\]
Observe that the matrix $ST^{-{1}/{2}}-S^*T^{{1}/{2}}$ (where $M^*$ denotes the
transpose of a matrix $M$) augments to an invertible
matrix. We call it the {\it Alexander matrix} corresponding to 
the Seifert matrix $S$, and denote it by $\Lambda_S$.

\begin{thm}\label{abelianthm}{\rm\cite{GK}}\qua
For $L$ a 0-framed 
boundary link in $M$, a $\ZHS$, and $S$ a Seifert matrix
for $L$, the value of
$\ZLMO(M,L)$, an element of the space $\JSS$, is equal to
\begin{equation}
\left(\nu(k_1)\sqcup \ldots \sqcup \nu(k_\mu)\right) \sqcup 
\Psi(
\Lambda_S)
^{-\frac{1}{2}}\sqcup(1+ \ldots\ ),\label{boundarywheels}
\end{equation}
where the error in the bracket above 
is a series of diagrams in $\JS$ corresponding to 
a series of generating diagrams each connected 
component of which has
at least two loops.
\end{thm}

\begin{rems}

(1)\qua That is, the error is a series of terms like
$\setlength{\unitlength}{0.75cm}
\thetdiagmz$ (supressing labels), i.e. from closed trivalent graphs
with at least 2 loops.

(2)\qua Let us recall how and where this is proved in the literature.
Observe that given a collection of Seifert surfaces for $L$, the
factor $\Psi(\Lambda_S)$ does not depend on the choice of a basis of
curves. It suffices, then, to prove the theorem for some special
choice of curves. The special choice of curves the proof uses 
is made by putting the 
surface in some standard disc-band form,
and taking the associated standard collection of curves. 
Now, the relevant piece of the LMO invariant is expressed
in terms of $\Psi$ of the equivariant linking matrix of a
('nice') surgery presentation. Here, one chooses
a certain surgery presentation that is canonically associated to a
disc-band decomposed Seifert surface: the so-called Y-view of
boundary links (see e.g. \cite{GGP}). It turns out
that $\Psi$ of the equivariant
linking matrix of this surgery presentation
can then be manipulated
to produce the above function of the Seifert matrix corresponding
to the chosen basis of curves. 
(Elements of this strategy were introduced in \cite{KS}.)

(3)\qua $\Psi(\Lambda_S)$ {\it does} depend on the isotopy
class of the collection. Its image in the space $\JSS$, however,
does not. These issues, and the above proof, 
are discussed in full detail in $\cite{GK}$. 
We remark that, in the form presented, 
this theorem does not depend on the two pieces of heavy machinery
employed by \cite{GK} -- it depends on neither
the adapted Kirby-Fenn-Rourke theorem nor the \cite{BLT} calculation
of the Kontsevich integral of the unknot.

(4)\qua The notation $\nu(k)$ denotes the Kontsevich integral of the unknot in 
${\cal A}(\star_{k})$ (see the Appendix). 
The space $\JSS$ is the quotient of $\JS$ obtained by
imposing the notorious
{\it link relations} (see \cite{BGRT}). Note that, for starters,
the $\sqcup$-product
is not well-defined in the presence of such relations. 
The meaning of the above equation, then,
is the following: the factors lie in $\JS$, by definition, and they are
to be multiplied in that space and then the result is to be projected
to $\JSS$.\end{rems}

\section{Casson's invariant}

We now recall the identification of Casson's formula with the degree
1 piece of the LMO invariant. We are 
going to do this in a way that is far from the easiest, as a warm-up and 
to provide context for the proof of our main theorem.

So, we wish to calculate the effect on the degree 1 piece of the LMO
invariant of a surgery on a knot $K$ in an integral homology
sphere $M$. In this case,
the determinant-like properties
give
\[
\Psi(ST^{-1/2}-S^*T^{1/2}) = \Psi(A_K(e^{k})).
\]
and we may rewrite equation \ref{boundarywheels}:
\begin{cor}
For a pair of an
$f$-framed knot $K$ in $M$, 
$\nu(k) \# \ZLMO(M,K)$ is equal to
\begin{equation}\label{leading}
\mbox{exp}_{\sqcup}\left( \frac{f}{2} 
\begin{picture}(1,0.7)(-0.25,0)
\qbezier[12](0,0)(0,0.5)(0.25,0.5)
\qbezier[12](0.25,0.5)(0.5,0.5)(0.5,0)
\put(-0.1,-0.4){$k$}
\put(0.4,-0.4){$k$}
\end{picture}
\right)\sqcup \left(1 + 
\frac{\lambda(M)}{2} \thetdiagsm +
(2\kappa -\frac{1}{2}
a_2(K))
\begin{picture}(1.2,0.5)(0.9,0.15)
\qbezier[6](1.25,0.75)(1.25,1)(1.5,1)
\qbezier[6](1.25,0.75)(1.25,0.5)(1.5,0.5)
\qbezier[6](1.75,0.75)(1.75,1)(1.5,1)
\qbezier[6](1.75,0.75)(1.75,0.5)(1.5,0.5)
\qbezier[6](1.25,0.75)(1,0.75)(1,0.5)
\qbezier[6](1,0.5)(1,0.375)(1.125,0.25)
\qbezier[6](1.125,0.25)(1.25,0.125)(1.25,0)
\qbezier[6](1.75,0.75)(2,0.75)(2,0.5)
\qbezier[6](2,0.5)(2,0.375)(1.875,0.25)
\qbezier[6](1.875,0.25)(1.75,0.125)(1.75,0)
\put(1.15,-0.4){$k$}
\put(1.65,-0.4){$k$}
\end{picture}
\, +\
\ldots \ \right),
\end{equation}
where the error inside the bracket is a series of diagrams
each of which has more than 2 trivalent vertices, and
$\kappa$ denotes the coefficient of the wheel with 2 legs
in $\nu(k)$.
\end{cor}

Two
multiplicative corrections were required 
to obtain this formula
from equation \ref{boundarywheels}: 
First, we had to multiply in a
copy $\nu(k)$ to allow use in 
the LMO surgery formula (equation \ref{LMOSURG}). 
Then,
we had to adjust the
framing of $K$ from 0 to $f$. 
See the Appendix for
some details on an elegant way to perform such corrections (using
the beautiful wheeling isomorphism).  

How to put this formula to use?
Let us briefly recall how the 
degree $\leq n$ part of $\ZLMO(M_{K})$ is determined 
from the series $\ZLMO(M,K)$. This uses
the map
\[
\iota_n : {\cal A}(\star_{k}) \rightarrow {\cal A}(\phi).
\]
We recall this by calculating:
\[
\iota_2 \left(
\begin{picture}(2.5,1.2)(-0.25,0)
\qbezier[12](0,0)(0,0.5)(0.25,0.5)
\qbezier[12](0.25,0.5)(0.5,0.5)(0.5,0)
\put(-0.1,-0.4){$k$}
\put(0.4,-0.4){$k$}
\put(1.15,-0.4){$k$}
\put(1.65,-0.4){$k$}
\qbezier[10](1.25,0.75)(1.25,1)(1.5,1)
\qbezier[10](1.25,0.75)(1.25,0.5)(1.5,0.5)
\qbezier[10](1.75,0.75)(1.75,1)(1.5,1)
\qbezier[10](1.75,0.75)(1.75,0.5)(1.5,0.5)
\qbezier[6](1.25,0.75)(1,0.75)(1,0.5)
\qbezier[6](1,0.5)(1,0.375)(1.125,0.25)
\qbezier[6](1.125,0.25)(1.25,0.125)(1.25,0)
\qbezier[6](1.75,0.75)(2,0.75)(2,0.5)
\qbezier[6](2,0.5)(2,0.375)(1.875,0.25)
\qbezier[6](1.875,0.25)(1.75,0.125)(1.75,0)
\end{picture}
\right).
\]
The map $\iota_n$
sends a diagram with other than $2n$ legs to zero. 
Otherwise, (our example, for example) 
there are 3 steps. First, glue $n$ chords into the legs, in all possible ways:
\[
\begin{picture}(2.5,1)(-0.25,-0)
\qbezier[12](0,0)(0,0.5)(0.25,0.5)
\qbezier[12](0.25,0.5)(0.5,0.5)(0.5,0)
\qbezier[10](1.25,0.75)(1.25,1)(1.5,1)
\qbezier[10](1.25,0.75)(1.25,0.5)(1.5,0.5)
\qbezier[10](1.75,0.75)(1.75,1)(1.5,1)
\qbezier[10](1.75,0.75)(1.75,0.5)(1.5,0.5)
\qbezier[6](1.25,0.75)(1,0.75)(1,0.5)
\qbezier[6](1,0.5)(1,0.375)(1.125,0.25)
\qbezier[6](1.125,0.25)(1.25,0.125)(1.25,0)
\qbezier[6](1.75,0.75)(2,0.75)(2,0.5)
\qbezier[6](2,0.5)(2,0.375)(1.875,0.25)
\qbezier[6](1.875,0.25)(1.75,0.125)(1.75,0)
\qbezier[12](0,0)(0,-0.5)(0.25,-0.5)
\qbezier[12](0.25,-0.5)(0.5,-0.5)(0.5,0)
\qbezier[12](1.25,0)(1.25,-0.5)(1.5,-0.5)
\qbezier[12](1.5,-0.5)(1.75,-0.5)(1.75,0)
\end{picture}
+
\begin{picture}(2.5,1)(-0.25,-0)
\qbezier[12](0,0)(0,0.5)(0.25,0.5)
\qbezier[12](0.25,0.5)(0.5,0.5)(0.5,0)
\qbezier[10](1.25,0.75)(1.25,1)(1.5,1)
\qbezier[10](1.25,0.75)(1.25,0.5)(1.5,0.5)
\qbezier[10](1.75,0.75)(1.75,1)(1.5,1)
\qbezier[10](1.75,0.75)(1.75,0.5)(1.5,0.5)
\qbezier[6](1.25,0.75)(1,0.75)(1,0.5)
\qbezier[6](1,0.5)(1,0.375)(1.125,0.25)
\qbezier[6](1.125,0.25)(1.25,0.125)(1.25,0)
\qbezier[6](1.75,0.75)(2,0.75)(2,0.5)
\qbezier[6](2,0.5)(2,0.375)(1.875,0.25)
\qbezier[6](1.875,0.25)(1.75,0.125)(1.75,0)
\qbezier[12](0,0)(0,-0.5)(0.6125,-0.5)
\qbezier[12](0.6125,-0.5)(1.25,-0.5)(1.25,0)
\qbezier[12](0.5,0)(0.5,-0.6)(1.1125,-0.6)
\qbezier[12](1.1125,-0.6)(1.75,-0.6)(1.75,0)
\end{picture}
+
\begin{picture}(2.5,1)(-0.25,-0)
\qbezier[12](0,0)(0,0.5)(0.25,0.5)
\qbezier[12](0.25,0.5)(0.5,0.5)(0.5,0)
\qbezier[10](1.25,0.75)(1.25,1)(1.5,1)
\qbezier[10](1.25,0.75)(1.25,0.5)(1.5,0.5)
\qbezier[10](1.75,0.75)(1.75,1)(1.5,1)
\qbezier[10](1.75,0.75)(1.75,0.5)(1.5,0.5)
\qbezier[6](1.25,0.75)(1,0.75)(1,0.5)
\qbezier[6](1,0.5)(1,0.375)(1.125,0.25)
\qbezier[6](1.125,0.25)(1.25,0.125)(1.25,0)
\qbezier[6](1.75,0.75)(2,0.75)(2,0.5)
\qbezier[6](2,0.5)(2,0.375)(1.875,0.25)
\qbezier[6](1.875,0.25)(1.75,0.125)(1.75,0)
\qbezier[12](0,0)(0,-0.75)(0.875,-0.75)
\qbezier[12](0.875,-0.75)(1.75,-0.75)(1.75,0)
\qbezier[12](0.5,0)(0.5,-0.4)(0.875,-0.4)
\qbezier[12](0.875,-0.4)(1.25,-0.4)(1.25,0)
\end{picture}.
\]
\vspace{0.1cm}

Then, replace any vertex-free loops
that have arisen with multiplicative factors of $-2n$, and, finally, set any
resulting diagrams with more than $2n$ vertices to zero. Thus, our 
example yields $-2\theta$. We will presently need the following result,
which is fun to prove \cite{LeG} ($D$ is any diagram of degree less than
or equal to $n$ with exactly 2 legs):
\begin{equation}\label{leform}
\iota_n
\left(
\begin{picture}(6,1)(-0.25,0.35)
\qbezier[12](0,0)(0,0.5)(0.25,0.5)
\qbezier[12](0.25,0.5)(0.5,0.5)(0.5,0)
\put(-0.1,-0.4){$k$}
\put(0.4,-0.4){$k$}
\put(2.9,-0.4){$k$}
\put(3.4,-0.4){$k$}
\qbezier[5](1,0.25)(1.75,0.25)(2.5,0.25)
\qbezier[12](3,0)(3,0.5)(3.25,0.5)
\qbezier[12](3.25,0.5)(3.5,0.5)(3.5,0)
\put(0.1,0.75){$\overbrace{\hspace{3.25cm}}^{n-1}$}
\put(4.5,0){$
\begin{picture}(1,1)(1.25,0)
\qbezier(1.25,0.75)(1.25,1)(1.5,1)
\qbezier(1.25,0.75)(1.25,0.5)(1.5,0.5)
\qbezier(1.75,0.75)(1.75,1)(1.5,1)
\qbezier(1.75,0.75)(1.75,0.5)(1.5,0.5)
\qbezier[6](1.25,0.75)(1,0.75)(1,0.5)
\qbezier[6](1,0.5)(1,0.375)(1.125,0.25)
\qbezier[6](1.125,0.25)(1.25,0.125)(1.25,0)
\qbezier[6](1.75,0.75)(2,0.75)(2,0.5)
\qbezier[6](2,0.5)(2,0.375)(1.875,0.25)
\qbezier[6](1.875,0.25)(1.75,0.125)(1.75,0)
\put(1.15,-0.4){$k$}
\put(1.65,-0.4){$k$}
\put(1.325,0.625){$D$}
\end{picture} 
$}
\end{picture}\right)=(-1)^{n-1}2^{n-1}(n-1)!\
\begin{picture}(1,1)(1,0.25)
\qbezier(1.25,0.75)(1.25,1)(1.5,1)
\qbezier(1.25,0.75)(1.25,0.5)(1.5,0.5)
\qbezier(1.75,0.75)(1.75,1)(1.5,1)
\qbezier(1.75,0.75)(1.75,0.5)(1.5,0.5)
\qbezier[6](1.25,0.75)(1,0.75)(1,0.5)
\qbezier[6](1,0.5)(1,0.375)(1.125,0.25)
\qbezier[6](1.125,0.25)(1.25,0.125)(1.25,0)
\qbezier[6](1.75,0.75)(2,0.75)(2,0.5)
\qbezier[6](2,0.5)(2,0.375)(1.875,0.25)
\qbezier[6](1.875,0.25)(1.75,0.125)(1.75,0)
\qbezier[4](1.25,-0.0)(1.25,-0.15)(1.5,-0.15)
\qbezier[4](1.5,-0.15)(1.75,-0.15)(1.75,-0.0)
\put(1.325,0.625){$D$}
\end{picture} 
\end{equation}

The degree $\leq n$ part of $\ZLMO(M_K)$ is given by the degree $\leq n$
part of the expression
\cite{LMO,LeD}:
\begin{equation}\label{LMOSURG}
\frac{
\iota_n\left( \ZLMO(M,K) \# \nu(k) \right)}
{ \iota_n\left( \ZLMO(S^3,U_f) \# \nu(k) \right) } \in {\cal A}(\phi),
\end{equation}
where $U_f$ is an unknot with the framing of $K$, and we have to take
the degree less than or equal to $n$ part of the result. 

We can now observe the effect of surgery on Casson's invariant.
Setting $n$ to 1 leads to a quick calculation, 
but we will learn more (on the
way to our main theorem) if we do this for arbitrary (positive) $n$.
So, substituting equation \ref{leading} into this expression, 
we find that $1 + \frac{\lambda(M_K)}{2}\theta + \ldots$ equals
\setlength{\unitlength}{.57cm}
\[
\frac{
\iota_n\left(
\frac{1}{n!}
\left(
\frac{f}{2} 
\begin{picture}(1,0.7)(-0.25,0)
\qbezier[12](0,0)(0,0.5)(0.25,0.5)
\qbezier[12](0.25,0.5)(0.5,0.5)(0.5,0)
\put(-0.1,-0.4){$k$}
\put(0.4,-0.4){$k$}
\end{picture}
\right)^{n}
\left(
1 + 
\frac{\lambda(M)}{2} \thetdiagsm \right) +
\frac{1}{(n-1)!}
\left(
\frac{f}{2} 
\begin{picture}(1,0.7)(-0.25,0)
\qbezier[12](0,0)(0,0.5)(0.25,0.5)
\qbezier[12](0.25,0.5)(0.5,0.5)(0.5,0)
\put(-0.1,-0.4){$k$}
\put(0.4,-0.4){$k$}
\end{picture}
\right)^{n-1}
(2\kappa -\frac{1}{2}
a_2(K))
\begin{picture}(1.2,1)(0.9,0.15)
\qbezier[6](1.25,0.75)(1.25,1)(1.5,1)
\qbezier[6](1.25,0.75)(1.25,0.5)(1.5,0.5)
\qbezier[6](1.75,0.75)(1.75,1)(1.5,1)
\qbezier[6](1.75,0.75)(1.75,0.5)(1.5,0.5)
\qbezier[6](1.25,0.75)(1,0.75)(1,0.5)
\qbezier[6](1,0.5)(1,0.375)(1.125,0.25)
\qbezier[6](1.125,0.25)(1.25,0.125)(1.25,0)
\qbezier[6](1.75,0.75)(2,0.75)(2,0.5)
\qbezier[6](2,0.5)(2,0.375)(1.875,0.25)
\qbezier[6](1.875,0.25)(1.75,0.125)(1.75,0)
\put(1.15,-0.4){$k$}
\put(1.65,-0.4){$k$}
\end{picture}
+\ \ldots\ 
\right) }
{
\iota_n\left(
\frac{1}{n!}
\left(
\frac{f}{2} 
\begin{picture}(1,0.7)(-0.25,0)
\qbezier[12](0,0)(0,0.5)(0.25,0.5)
\qbezier[12](0.25,0.5)(0.5,0.5)(0.5,0)
\put(-0.1,-0.4){$k$}
\put(0.4,-0.4){$k$}
\end{picture}
\right)^{n}
+ 
\frac{1}{(n-1)!}
\left(
\frac{f}{2} 
\begin{picture}(1,0.7)(-0.25,0)
\qbezier[12](0,0)(0,0.5)(0.25,0.5)
\qbezier[12](0.25,0.5)(0.5,0.5)(0.5,0)
\put(-0.1,-0.4){$k$}
\put(0.4,-0.4){$k$}
\end{picture}
\right)^{n-1}
2\kappa
\begin{picture}(1.2,1)(0.9,0.15)
\qbezier[6](1.25,0.75)(1.25,1)(1.5,1)
\qbezier[6](1.25,0.75)(1.25,0.5)(1.5,0.5)
\qbezier[6](1.75,0.75)(1.75,1)(1.5,1)
\qbezier[6](1.75,0.75)(1.75,0.5)(1.5,0.5)
\qbezier[6](1.25,0.75)(1,0.75)(1,0.5)
\qbezier[6](1,0.5)(1,0.375)(1.125,0.25)
\qbezier[6](1.125,0.25)(1.25,0.125)(1.25,0)
\qbezier[6](1.75,0.75)(2,0.75)(2,0.5)
\qbezier[6](2,0.5)(2,0.375)(1.875,0.25)
\qbezier[6](1.875,0.25)(1.75,0.125)(1.75,0)
\put(1.15,-0.4){$k$}
\put(1.65,-0.4){$k$}
\end{picture}
\ +\ \ldots \right)\hspace{2cm} }.
\setlength{\unitlength}{1cm}
\]
Here is what is happening: 
terms which do not have precisely $2n$ legs are killed by
$\iota_n$. 
Assuming 
$\iota_n$ is to be operated on an expression 
of the form exp$_\sqcup(\mbox{chords})\sqcup X$, 
for $X$ some series of diagrams,
then to each term in $X$ we multiply as many chords as we need to get the
right number of legs, and that is what we see above. The error terms 
above 
arise from diagrams in $X$ with more than 2 trivalent vertices:
acting with $\iota_n$ on such a diagram (producted with an appropriate
number of chords)
inevitably yields a diagram with more than 2 loops. 
Evaluating returns the expected 
$1 + \left( \lambda(M) + f a_2(M,K) \right) \frac{1}{2}\theta + \mbox{diagrams
with $>2 $ loops}$\ .

\section{The loop expansion and rationality}

So, the 1-loop part of $\ZLMO(M,K)$ is a polynomial, more or less. This
and general
principles let us write $\ZLMO(M,K)$ as
\[
\underbrace{
\nu(k)\sqcup 
\left(\Psi(A_K(e^k))\right)^{-\frac{1}{2}}
\hspace{-0.4cm}
\begin{array}{l}
\\ \\
\end{array}
}_{
\mbox{The wheels}}
\sqcup\,
\mbox{exp}_{\sqcup}
(
\underbrace{
\
R_2 + R_3 + R_4 + \ldots\ 
\hspace{-0.4cm}
\begin{array}{l}
\\ \\
\end{array}
}_{\mbox{$>1$ loops}}),
\]
where $R_l$ is a (unique modulo IHX and AS)
series of connected diagrams each with $l$ loops. Lev
Rozansky's studies of the Melvin-Morton expansion of the coloured
Jones polynomial \cite{RMM}, combined with results in the theory of
finite-type invariants (for example the proof of Melvin-Morton \cite{BG}), 
led
him to conjecture that each piece $R_l$ might be generated from
``finite'' data as well: from a finite collection of
1-variable polynomials \cite{RozG}. 
\begin{thm}\label{rationality}{\rm\cite{K}}\qua
There exists a finite $\Qset$-linear
combination of generating diagrams with $l$ loops, 
each of whose edges is labelled
by exactly one power series of the form
\[
\frac{\mbox{\rm Laurent polynomial in $e^k$}}{A_K(e^{k})},
\]
which generates the series $R_l$.
\end{thm}
\paragraph{The 2-loop piece}
Let us focus now on the 2-loop piece, $R_2$. According
to this theorem, and the following exercise (for rational functions
$\{q_i(t)\}$ non-singular at 1),
\setlength{\unitlength}{0.5cm}
\begin{eqnarray*}
& & 
\begin{picture}(10,1)(-2,0.35)
\qbezier[8](0,0)(0,1)(1,1)
\qbezier[8](1,1)(2,1)(2,0)
\qbezier[8](0,0)(0,-1)(1,-1)
\qbezier[8](1,-1)(2,-1)(2,0)
\qbezier[8](4,0)(4,1)(5,1)
\qbezier[8](5,1)(6,1)(6,0)
\qbezier[8](4,0)(4,-1)(5,-1)
\qbezier[8](5,-1)(6,-1)(6,0)
\qbezier[12](2,0)(3,0)(4,0)
\put(0,0){\circle*{0.2}}
\put(0,-0.3){\vector(0,1){0.7}}
\put(-2.3,-0.1){$q_1(e^k)$}
\put(2.05,0.7){$q_2(e^k)$}
\put(6.25,-0.1){$q_3(e^k)$}
\put(6,-0.3){\vector(0,1){0.7}}
\put(2.7,0){\vector(1,0){0.7}}
\put(3,0){\circle*{0.2}}
\put(6,0){\circle*{0.2}}
\end{picture} \\
& & \\
& &
\begin{picture}(18,1)(-11,0.35)
\qbezier[8](0,0)(0,1)(1,1)
\qbezier[8](1,1)(2,1)(2,0)
\qbezier[8](0,0)(0,-1)(1,-1)
\qbezier[8](1,-1)(2,-1)(2,0)
\qbezier[14](1,1)(1,0)(1,-1)
\put(0,0){\circle*{0.2}}
\put(0,-0.3){\vector(0,1){0.7}}
\put(2,0){\circle*{0.2}}
\put(2,-0.3){\vector(0,1){0.7}}
\put(2.3,-0.1){$(q_3(e^k)-q_3(e^{-k}))\ ,\ $}
\put(-10.5,-0.1){$=\ \ \ q_2(1) \times (q_1(e^k)-q_1(e^{-k}))$}
\end{picture}
\end{eqnarray*}
\setlength{\unitlength}{0.8cm}

there exists a finite set of 4-tuples  
$\{\lambda_i, p_i(t), q_i(t), r_i(t)\}$
of a rational $\lambda_i$ together with three polynomials, such that
\begin{equation}
\begin{picture}(7,1.5)(-3.75,-0.5)
\put(-4.5,0){$R_2 = \sum_i \lambda_i$}
\qbezier[16](0,0)(0,1)(1,1)
\qbezier[16](3,0)(3,1)(2,1)
\qbezier[16](0,0)(0,-1)(1,-1)
\qbezier[16](3,0)(3,-1)(2,-1)
\qbezier[20](1,-1)(1,0)(1,1)
\qbezier[12](1,1)(1.5,1)(2,1)
\qbezier[12](1,-1)(1.5,-1)(2,-1)
\put(0,0){\circle*{0.15}}
\put(1,0){\circle*{0.15}}
\put(3,0){\circle*{0.15}}
\put(0,-0.2){\vector(0,1){0.5}}
\put(1,-0.2){\vector(0,1){0.5}}
\put(3,-0.2){\vector(0,1){0.5}}
\put(-1.75,-0){$
\frac{p_i(e^k)}{A_K(e^k)}$}
\put(1.025,-0){$
\frac{q_i(e^k)}{A_K(e^k)}$}
\put(3.025,-0){$
\frac{r_i(e^k)}{A_K(e^k)}$}
\put(5,0){$.$}
\end{picture}
\begin{array}{l} \\ \\ \\
\end{array}
\end{equation}

We digress for 2 paragraphs on the the question of how to present 
this result.  
It is possible, at this point, to pass to a space of trivalent
graphs with edges labelled by rational functions of $t$
(that is, with no mention of power series). For example,
in \cite{GK} the authors defined a space, ${\cal A}(\Lambda_{loc})$,
which is generated by such diagrams, and which has certain relations
which allow one to ``Push'' around the labels.

Keeping in mind our aim of presenting this work ``in the simplest
possible terms'' we will not discuss these issues here
(see \cite{GK} for a full discussion). That is, we will
continue to present results at the level of declaring that
some element in question is an element of $\TS$
generated by some finite combination of trivalent diagrams with 
edges labelled by power series in $k$ generated by rational
functions in $e^k$. Suffice it to say that, at least for the 2 loop
piece, such an expression uniquely determines some polynomials (up to
some symmetries)%
\footnote{Note added: this is not true for a general number of loops.
The question of injectivity has recently been resolved
\cite{BPM} (these spaces do not inject into $\TS$).}.

Returning to the main discussion:
Let us write $R_2$, the 2-loop invariant, as 
$\Theta(M,K) \in {\cal A}(\star_k).$ This invariant generalises 
Casson's invariant in the sense that, obviously,
\begin{equation}\label{speceqn}
\left.
\left(
\begin{array}{l} \\
\end{array}
\hspace{-0.3cm}
\Theta(M,K)\,
\right)
\right|_{k=0}
=
\frac{\lambda(M)}{2} \thetdiagsm\ ,
\end{equation}
in other words, when diagrams with legs are set to zero. That this 
generalisation is natural is well illustrated by our main
theorem:

\section{Surgery on a sublink of a boundary link}

We come, at last, to our generalisation to $\Theta(M,K)$ of Casson's
formula (equation \ref{casson}). Ideally, we would like to describe
the effect on $\Theta(M,K)$ of a surgery on {\it any} $\pm 1$-framed
knot $K'$ in $M-K$ (this, for example, would let us change crossings).
Here, however, we can only give a formula describing the effect of
a surgery on a $\pm 1$-framed 
knot $K'$ with the property that 
$(K,K')$ forms a {\it boundary link}
in $M$. The pair that results we denote $(M,(K,K'))_{K'}$.

There are several, related, reasons why this is a natural class of surgeries
to consider. Firstly, to
anticipate the proof: if $(K,K')$ forms a boundary link, then the invariant
$\ZLMO(M,(K,K'))$ can be written {\it without trees}, which lets
us control
the contributions to $\Theta(M,K)$,
the 2-loop piece, when we evaluate the LMO surgery formula. 
Secondly, recall 
that $\Theta(M,K)$ appears in $\ZLMO(M,K)$ as
\[
\nu(k)\sqcup 
\left(\Psi(A_K(e^k))\right)^{-\frac{1}{2}}\sqcup
\left(1 + \Theta(M,K) + \ldots \right).
\]
A general principle of quantum topology then suggests that it should be
possible to compare $\Theta(M_1,K_1)$ and $\Theta(M_2,K_2)$
when their Alexander polynomials coincide. As it happens, surgery
on a $\pm 1$-framed component $K'$ of a boundary link $(K,K')$ does
not affect the Alexander polynomial of $K$ (this is because 
the lift of $K'$ bounds
in the universal cyclic cover of $K$, for example).

The statement uses the following notation, where $D$ is a
diagram in ${\cal A}(\star_{\{k,k'\}})$ (possibly with some number
of $k$-labelled legs, not shown). 
\setlength{\unitlength}{1.3cm}
\begin{equation}
\left< 
\frac{1}{2}
\begin{picture}(1,0.7)(-0.25,-0.125)
\qbezier[16](0,-0.1)(0,0.25)(0.25,0.25)
\qbezier[16](0.25,0.25)(0.5,0.25)(0.5,-0.1)
\put(-0.1,-0.4){$k'$}
\put(0.4,-0.4){$k'$}
\end{picture} ,
\begin{picture}(1.5,0.7)(-0.25,0)
\qbezier(0,0.5)(0,0.75)(0.25,0.75)
\qbezier(0,0.5)(0,0.25)(0.25,0.25)
\qbezier(1,0.5)(1,0.75)(0.75,0.75)
\qbezier(1,0.5)(1,0.25)(0.75,0.25)
\qbezier(0.25,0.75)(0.5,0.75)(0.75,0.75)
\qbezier(0.25,0.25)(0.5,0.25)(0.75,0.25)
\qbezier[12](0.25,0.25)(0.25,0.125)(0.25,-0.1)
\qbezier[12](0.75,0.25)(0.75,0.125)(0.75,-0.1)
\qbezier[3](0.35,0)(0.5,0)(0.65,0)
\put(0.15,-0.4){$k'$}
\put(0.65,-0.4){$k'$}
\put(0.35,0.4){$D$}
\put(0.125,-0.55){$\underbrace{\hspace{1cm}}_n$}
\end{picture}
\right>
=
\left\{
\begin{array}{ll}
\begin{picture}(1.5,0.7)(0,0)
\qbezier(0,0.5)(0,0.75)(0.25,0.75)
\qbezier(0,0.5)(0,0.25)(0.25,0.25)
\qbezier(1,0.5)(1,0.75)(0.75,0.75)
\qbezier(1,0.5)(1,0.25)(0.75,0.25)
\qbezier(0.25,0.75)(0.5,0.75)(0.75,0.75)
\qbezier(0.25,0.25)(0.5,0.25)(0.75,0.25)
\qbezier[8](0.25,0.25)(0.25,0.125)(0.25,0)
\qbezier[8](0.75,0.25)(0.75,0.125)(0.75,0)
\put(0.35,0.4){$D$}
\qbezier[10](0.25,-0.)(0.25,-0.3)(0.5,-0.3)
\qbezier[10](0.75,-0.)(0.75,-0.3)(0.5,-0.3)
\end{picture}
&
\mbox{if $n=2$,} \\
& \\
\hspace{0.5cm} 0 & \mbox{otherwise.}
\end{array} 
\right.
\end{equation}

Instead of taking the coefficient of $k^2$ in a series
associated to the Alexander polynomial (as appears in Casson's
formula), we use this operation to ``take the coefficient
of $(k')^2$'' in the diagram-valued determinant associated to the boundary 
link $(K,K')$. To be precise:

\begin{thm}[The Main Theorem]\label{mainth}

\

Let $(M,(K,K'))$ be a $(0,f)$-framed boundary link $(K,K')$ in 
$M$, a $\ZHS$, (where $f$ is plus or minus 1,)
and let $S$ be a Seifert matrix corresponding
to some choice of a disjoint pair of
Seifert surfaces $(\Sigma,\Sigma')$ for $(K,K')$. Then,
\begin{equation}\label{maineqn}
\Theta\left((M,(K,K'))_{K'}\right) =
\Theta\left( M,K \right) 
-\ f\,
\left< 
\frac{1}{2}
\begin{picture}(1,0.7)(-0.25,-0.125)
\qbezier[16](0,-0.1)(0,0.25)(0.25,0.25)
\qbezier[16](0.25,0.25)(0.5,0.25)(0.5,-0.1)
\put(-0.1,-0.4){$k'$}
\put(0.4,-0.4){$k'$}
\end{picture} \ ,\
\Psi
\left(
W
\right)^{-1/2}
\right>,
\end{equation}
where 
\[
W= \Lambda_S(k,k') \Lambda_S(k,0)^{-1}.
\]
\end{thm}
Recall that $\Lambda_S(k,k')$ 
denotes the Alexander matrix corresponding to $S$: 
\[ 
\Lambda_S(k,k')=
ST(k,k')^{-1/2} - S^*T(k,k')^{1/2}.\]
Three remarks:

(1)\qua Let $S_K$ be the Seifert matrix of $K$ by itself.
Observe that it follows from
the multiplicativity of the determinant that the exponential inside
the pairing can be written
\[
\left(\frac{\Psi\left(\Lambda_S\right)}
{\Psi\left(\Lambda_{S_K}\right)}\right)^{-\frac{1}{2}}.
\]

(2)\qua If we are to believe this formula, then it must, at least, give us
a ``rational'' $\Theta((M,(K,K'))_{K'})$ (see Theorem \ref{rationality}). 
To see that this is the case, write $\Psi$ as exp\,tr\,log,
and consider the argument of the exponential:
\begin{equation}\label{argument}
\frac{1}{2} \sum_{p=1}^{\infty}
\WT\left( 
\frac{\left((\Lambda_S(k,0)-\Lambda_S(k,k'))
\Lambda_S(k,0)^{-1}\right)^{p}}{p} \right).
\end{equation}
First, we perform 
some cancellations to write this in terms of integral powers
of $e^{k}$. Let $\hat{\Lambda}_S= \Lambda_S T(-k/2 , 0)$, that is, $ \hat{\Lambda}_S = S T(-k,-k'/2) - S^* T(0,k'/2)$.
Observe that the expression (\ref{argument}) may be
written:
\[
\frac{1}{2} \sum_{p=1}^{\infty}
\WT\left( 
\frac{\left((\hat{\Lambda}_S(k,0)-\hat{\Lambda}_S(k,k'))
\hat{\Lambda}_S(k,0)^{-1}\right)^{p}}{p} \right).
\]

Now, note that
the factor $
(\hat{\Lambda}_S(k,0)-\hat{\Lambda}_S(k,k'))$ is independent of $k$
and is divided by $k'$. 
In fact, 
for $M_a$ and $M_b$ matrices of integers:
\[
(\hat{\Lambda}_S(k,0)-\hat{\Lambda}_S(k,k')) = k'M_a + {k'}^2 M_b + (\mbox{terms with $\geq 3$ factors $k'$}).
\]
Consider now the second factor, $\hat{\Lambda}_S(k,0)^{-1}$.
Observe that this factor has the form (for $B$
and $C$ some matrices of integers, $|C|=1$, and $\hat{\Lambda}_{S_K}(k)$ denoting
$(S_K e^{-k} - S_K^*)$):
\[
\left[
\begin{array}{cc}
\hat{\Lambda}_{S_K}(k)&
0 \\
 (e^{-k}-1) B & C
\end{array}
\right]^{-1} 
 = \frac{1}{A_K(e^{k})} M_c(e^{k}),
\]
where $M_c(e^{k})$ is some matrix of Laurent polynomials in $e^k$.

Now, only terms with less than 2 legs labelled $k'$ will contribute
to formula (\ref{maineqn}). 
Thus, the contribution of the $p=1$ term is precisely the sum of:
\vspace{-0.25cm}
\begin{equation}\label{diag}
\begin{picture}(1,0.75)(1.8,-0.1)
\put(-0.1,-0.05){$k'$}
\qbezier[11](0.3,0)(0.5,0)(0.7,0)
\qbezier[15](0.7,0)(0.7,0.4)(1,0.4)
\qbezier[15](1.3,0)(1.3,0.4)(1,0.4)
\qbezier[15](1.3,0)(1.3,-0.4)(1,-0.4)
\qbezier[15](0.7,0)(0.7,-0.4)(1,-0.4)
\put(1.3,0){\circle*{0.1}}
\put(1.3,-0.15){\vector(0,1){0.4}}
\put(1.5,-0.05){$\sum_{i,j}(M_a)_{ij} \frac{1}{A_K(e^{k})} 
(M_c(e^{k}))_{ji}$}
\end{picture}
\end{equation}
\vspace{-0.25cm}
and
\vspace{-0.25cm}
\begin{equation}\label{diaga}
\begin{picture}(1,0.75)(1.8,-0.1)
\put(-0.1,-0.25){$k'$}
\qbezier[11](0.3,-0.2)(0.5,-0.2)(0.7,-0.2)
\put(-0.1,0.15){$k'$}
\qbezier[11](0.3,0.2)(0.5,0.2)(0.7,0.2)
\qbezier[15](0.7,0)(0.7,0.4)(1,0.4)
\qbezier[15](1.3,0)(1.3,0.4)(1,0.4)
\qbezier[15](1.3,0)(1.3,-0.4)(1,-0.4)
\qbezier[15](0.7,0)(0.7,-0.4)(1,-0.4)
\put(1.3,0){\circle*{0.1}}
\put(1.3,-0.15){\vector(0,1){0.4}}
\put(1.5,-0.05){$\sum_{i,j}(M_b)_{ij} \frac{1}{A_K(e^{k})} 
(M_c(e^{k}))_{ji}$}
\end{picture}
\end{equation}
\vspace{-0.25cm}

For the purposes of the discussion, this is a good point
to make precise some terminology:
\begin{defn}\label{fragment}
A {\it fragment of a generating diagram} 
is a unitrivalent diagram with univalent vertices
labelled by $k'$, 
with edges possibly labelled by power series
in $k$, and such that each connected component has at least 1 trivalent 
vertex. A {\it fragment of a rational generating diagram} satisfies,
in addition, that the edge labels are rational functions
in $e^k$.
\end{defn}
Diagrams (\ref{diag}) and (\ref{diaga}) are 
clearly such things -- each a fragment of a rational
generating diagram.
 
The pairing $<,>$
will map a product of such diagrams to a rational generating diagram.
It is easy to see that the $p=2$ term is a sum of such fragments, 
and that terms
for $p\geq 3$ will all have at least $3$ legs labelled $k'$, 
and hence will be killed by the pairing $<,>$. 

(3)\qua Our theorem suggests studying the following ``finite-type'' filtration.
Let ${\bf MK}$ be the Abelian group freely generated by pairs $(M,K)$
of a knot $K$ in a $\ZHS$ $M$. Let ${\bf MK}_n^\partial$ denote the subgroup
generated by elements of ${\bf MK}$ corresponding
to pairs $(M,L)$ of a boundary link $L$ with a distinguished
component $K$ in $M$, a $\ZHS$: the correspondence is to take the
obvious
alternating sum $\sum_{L'\subset L-K}(-1)^\epsilon(M,(L',K))_{L'}$.
This filtration is clearly as steep as the ``loop'' filtration (see
\cite{GR}). 

\proof[Proof of the Main Theorem]

We now turn to the proof of the main theorem. This will be a step-by-step
copy of the proof of Casson's formula. First, in analogy
with Equation \ref{leading}, we need
an expression for $\ZLMO(M,(K.K'))$ that identifies the
contributing terms.

We start with Theorem \ref{abelianthm}. Let $K'_0$ be $K'$ with the
zero framing. Theorem \ref{abelianthm}, together with the fact that
\[ 
\ZLMO(M,(K,K'))|_{k'=0}=\ZLMO(M,K),
\]
tells us that 
$\ZLMO(M,(K,K'_0))$ is equal to:
\setlength{\unitlength}{1cm}
\begin{eqnarray*}
& & 
\hspace{-1.5cm} 
\nu(k) \sqcup 
\left(\Psi(\Lambda_S)\right)^{-\frac{1}{2}} 
\sqcup
\left( 1+\Theta(M,K) + 
\kappa 
\begin{picture}(1.2,0.8)(0.9,0.15)
\qbezier[6](1.25,0.75)(1.25,1)(1.5,1)
\qbezier[6](1.25,0.75)(1.25,0.5)(1.5,0.5)
\qbezier[6](1.75,0.75)(1.75,1)(1.5,1)
\qbezier[6](1.75,0.75)(1.75,0.5)(1.5,0.5)
\qbezier[6](1.25,0.75)(1,0.75)(1,0.5)
\qbezier[6](1,0.5)(1,0.375)(1.125,0.25)
\qbezier[6](1.125,0.25)(1.25,0.125)(1.25,0)
\qbezier[6](1.75,0.75)(2,0.75)(2,0.5)
\qbezier[6](2,0.5)(2,0.375)(1.875,0.25)
\qbezier[6](1.875,0.25)(1.75,0.125)(1.75,0)
\put(1.15,-0.4){$k'$}
\put(1.65,-0.4){$k'$}
\end{picture}
+\ \ldots\ \right),
\end{eqnarray*}
where the error term
is a series of fragments of generating diagrams
(see Defintion \ref{fragment})
with more than 2 trivalent vertices%
\footnote{We remark on a possible source of confusion here. This sentence
refers to the number of trivalent vertices in the generating diagram,
not to the number of trivalent vertices in the diagrams it generates.
For example: the fragment in diagram \ref{frag} has 4 trivalent vertices
while the series it generates contains diagrams with an arbitrarily high 
number of trivalent vertices (and $k$-labelled legs).}. 
Note that the maps $\{\iota_n\}$ send products of such terms 
with chords to
generating diagrams with {\it at least} 3 loops.
Let us examine an example. Under $\iota_2$, the following
product of a fragment of a generating diagram and a chord
is mapped:
\begin{equation}\label{frag}
\begin{picture}(1,0.7)(0.25,-0.8)
\qbezier[12](0,0)(0,0.5)(0.25,0.5)
\qbezier[12](0.25,0.5)(0.5,0.5)(0.5,0)
\put(-0.1,-0.4){$k'$}
\put(0.4,-0.4){$k'$}
\end{picture}
\begin{picture}(1.8,2.2)(0,-0.85)
\qbezier[12](0,0.5)(0,0.2)(0.1,0.1)
\qbezier[12](0.5,0)(0.2,0)(0.1,0.1)
\qbezier[12](0,0.5)(0,0.8)(0.1,0.9)
\qbezier[12](0.5,1)(0.2,1)(0.1,0.9)
\qbezier[12](1,0.5)(1,0.2)(0.9,0.1)
\qbezier[12](0.5,0)(0.8,0)(0.9,0.1)
\qbezier[12](1,0.5)(1,0.8)(0.9,0.9)
\qbezier[12](0.5,1)(0.8,1)(0.9,0.9)
\qbezier[12](0,0.5)(0.5,0.5)(1,0.5)
\qbezier[12](0.1,0.1)(-0.1,-0.1)(-0.1,-0.8)
\qbezier[12](0.9,0.1)(1.1,-0.1)(1.1,-0.8)
\put(1.1,0.7){\vector(-1,1){0.4}}
\put(0.9,0.9){\circle*{0.1}}
\put(-0.2,-1.2){$k'$}
\put(1.0,-1.2){$k'$}
\put(1,1){$s(k)$}
\end{picture}
\begin{picture}(1.9,0.5)(0,-0.8)
\put(0.25,0){$\mapsto\ -2\ \hspace{3cm} .$}
\end{picture}
\begin{picture}(1,2)(0,-0.85)
\qbezier[12](0,0.5)(0,0.2)(0.1,0.1)
\qbezier[12](0.5,0)(0.2,0)(0.1,0.1)
\qbezier[12](0,0.5)(0,0.8)(0.1,0.9)
\qbezier[12](0.5,1)(0.2,1)(0.1,0.9)
\qbezier[12](1,0.5)(1,0.2)(0.9,0.1)
\qbezier[12](0.5,0)(0.8,0)(0.9,0.1)
\qbezier[12](1,0.5)(1,0.8)(0.9,0.9)
\qbezier[12](0.5,1)(0.8,1)(0.9,0.9)
\qbezier[12](0,0.5)(0.5,0.5)(1,0.5)
\qbezier[12](0.1,0.1)(-0.1,-0.1)(-0.1,-0.4)
\qbezier[12](0.9,0.1)(1.1,-0.1)(1.1,-0.4)
\qbezier[12](-0.1,-0.4)(-0.1,-0.8)(0.5,-0.8)
\qbezier[12](1.1,-0.4)(1.1,-0.8)(0.5,-0.8)
\put(1.1,0.7){\vector(-1,1){0.4}}
\put(0.9,0.9){\circle*{0.1}}
\put(1,1){$s(k)$}
\end{picture}
\end{equation}

It is to be understood 
that the error terms in the equations that follow are also of
this form, and consequently do not contribute to the 2-loop piece, for the
same reason.

Following the Appendix, we can adjust 
the framing of $K'_0$ to $f$, and $\#$-multiply by a copy of $\nu(k')$ to
find that $\ZLMO(M,(K,K'))\#_{k'}(\nu(k'))$ is equal to 
\[
\mbox{exp}_{\sqcup}\left( \frac{f}{2} 
\begin{picture}(1,0.7)(-0.25,0)
\qbezier[12](0,0)(0,0.5)(0.25,0.5)
\qbezier[12](0.25,0.5)(0.5,0.5)(0.5,0)
\put(-0.1,-0.4){$k'$}
\put(0.4,-0.4){$k'$}
\end{picture}\,
\right)
\sqcup \nu(k) \sqcup 
\left(\Psi(\Lambda_S)\right)^{-\frac{1}{2}} 
\sqcup
\left(1
+\Theta(M,K)+
2\kappa 
\begin{picture}(1.2,0.5)(0.9,0.15)
\qbezier[6](1.25,0.75)(1.25,1)(1.5,1)
\qbezier[6](1.25,0.75)(1.25,0.5)(1.5,0.5)
\qbezier[6](1.75,0.75)(1.75,1)(1.5,1)
\qbezier[6](1.75,0.75)(1.75,0.5)(1.5,0.5)
\qbezier[6](1.25,0.75)(1,0.75)(1,0.5)
\qbezier[6](1,0.5)(1,0.375)(1.125,0.25)
\qbezier[6](1.125,0.25)(1.25,0.125)(1.25,0)
\qbezier[6](1.75,0.75)(2,0.75)(2,0.5)
\qbezier[6](2,0.5)(2,0.375)(1.875,0.25)
\qbezier[6](1.875,0.25)(1.75,0.125)(1.75,0)
\put(1.15,-0.4){$k'$}
\put(1.65,-0.4){$k'$}
\end{picture} + \ldots \right).
\]
This can be re-organised to isolate the factors with no legs labelled $k'$
and with no more than 2 loops:
\begin{eqnarray}
& & 
\hspace{-0.5cm}=
\left(\nu(k)\sqcup \Psi(\Lambda_{S_K})^{-\frac{1}{2}} 
\sqcup (1 + \Theta(M,K))\right) 
\nonumber \\
& & \hspace{0.5cm} \sqcup\, \mbox{exp}_{\sqcup}\left( \frac{f}{2} 
\begin{picture}(1,0.7)(-0.25,0)
\qbezier[12](0,0)(0,0.5)(0.25,0.5)
\qbezier[12](0.25,0.5)(0.5,0.5)(0.5,0)
\put(-0.1,-0.4){$k'$}
\put(0.4,-0.4){$k'$}
\end{picture}\,
\right)
\sqcup \left(
\frac{\Psi(\Lambda_S)}{\Psi(\Lambda_{S_K})} \right)^{-\frac{1}{2}}
\sqcup \left(1 + 2\kappa 
\begin{picture}(1.2,0.5)(0.9,0.15)
\qbezier[6](1.25,0.75)(1.25,1)(1.5,1)
\qbezier[6](1.25,0.75)(1.25,0.5)(1.5,0.5)
\qbezier[6](1.75,0.75)(1.75,1)(1.5,1)
\qbezier[6](1.75,0.75)(1.75,0.5)(1.5,0.5)
\qbezier[6](1.25,0.75)(1,0.75)(1,0.5)
\qbezier[6](1,0.5)(1,0.375)(1.125,0.25)
\qbezier[6](1.125,0.25)(1.25,0.125)(1.25,0)
\qbezier[6](1.75,0.75)(2,0.75)(2,0.5)
\qbezier[6](2,0.5)(2,0.375)(1.875,0.25)
\qbezier[6](1.875,0.25)(1.75,0.125)(1.75,0)
\put(1.15,-0.4){$k'$}
\put(1.65,-0.4){$k'$}
\end{picture} + \ldots \right). \label{justhere}
\end{eqnarray}

We will presently substitute this expression into the following formula,
which calculates, for some $n$,
the degree $\leq n$ part of $\ZLMO((M,(K,K'))_{K'})$
(see \cite{LMO,LeD,MO}):
\begin{equation}\label{calc}
\frac{\iota_n\left( \ZLMO(M,(K,K')) \#_{k'} \nu(k') \right)}
{\iota_n\left( \ZLMO(S^3,U_f) \#_k \nu(k) \right)}\ \in {\cal A}(\star_{k}),
\end{equation}
where the $\iota_n$ in the numerator acts on the labels $k'$. We will
observe that the 2-loop part of the result, 
as $n$ varies, is the degree $\leq n$ part of
some fixed series, which we can conclude calculates 
the 2-loop part of $\ZLMO((M,(K,K'))_{K'})$. 

Let us fix some $n$, then, and consider the numerator.
Denote the leading bracket of expression (\ref{justhere}) by $\alpha$.
Let $\psi$ denote the piece of
$\left(
\Psi(\Lambda_S)/\Psi(\Lambda_{S_K}) \right)^{-\frac{1}{2}}$ with precisely 2 legs
labelled by $k'$.
The numerator, then, is clearly the degree $\leq n$ piece of:
\setlength{\unitlength}{0.72cm}
\[
\iota_n\left(
\alpha\sqcup \left(
\frac{1}{n!}
\left(
\frac{f}{2} 
\begin{picture}(1,0.7)(-0.25,0)
\qbezier[12](0,0)(0,0.5)(0.25,0.5)
\qbezier[12](0.25,0.5)(0.5,0.5)(0.5,0)
\put(-0.1,-0.4){$k'$}
\put(0.4,-0.4){$k'$}
\end{picture}
\right)^{n}
+
\frac{1}{(n-1)!}
\left(
\frac{f}{2} 
\begin{picture}(1,0.7)(-0.25,0)
\qbezier[12](0,0)(0,0.5)(0.25,0.5)
\qbezier[12](0.25,0.5)(0.5,0.5)(0.5,0)
\put(-0.1,-0.4){$k'$}
\put(0.4,-0.4){$k'$}
\end{picture}
\right)^{n-1}
\sqcup
\left(\psi + 
2\kappa 
\begin{picture}(1.2,1)(0.9,0.15)
\qbezier[6](1.25,0.75)(1.25,1)(1.5,1)
\qbezier[6](1.25,0.75)(1.25,0.5)(1.5,0.5)
\qbezier[6](1.75,0.75)(1.75,1)(1.5,1)
\qbezier[6](1.75,0.75)(1.75,0.5)(1.5,0.5)
\qbezier[6](1.25,0.75)(1,0.75)(1,0.5)
\qbezier[6](1,0.5)(1,0.375)(1.125,0.25)
\qbezier[6](1.125,0.25)(1.25,0.125)(1.25,0)
\qbezier[6](1.75,0.75)(2,0.75)(2,0.5)
\qbezier[6](2,0.5)(2,0.375)(1.875,0.25)
\qbezier[6](1.875,0.25)(1.75,0.125)(1.75,0)
\put(1.15,-0.4){$k'$}
\put(1.65,-0.4){$k'$}
\end{picture}\right)
+\ \ldots\ 
\right) \right),
\]
\setlength{\unitlength}{1cm}
which evaluates (see equation \ref{leform}) 
to the degree $\leq n$ piece of
\[
(-1)^n f^n \alpha \sqcup \left( 1 - f^{-1} 
\left<\frac{1}{2}
\begin{picture}(1,0.7)(-0.25,-0.2)
\qbezier[16](0,-0.1)(0,0.25)(0.25,0.25)
\qbezier[16](0.25,0.25)(0.5,0.25)(0.5,-0.1)
\put(-0.1,-0.4){$k'$}
\put(0.4,-0.4){$k'$}
\end{picture} \ ,\
\psi
\right>
-2f^{-1}\kappa
\thetdiagsm
+ \ldots \right).
\]
The denominator, on the other hand, may be similarly calculated to be
the degree $\leq n$ part of:
\[
(-1)^n f^n \left( 1 - 2f^{-1}\kappa \thetdiagsm
+ \ldots \right).
\]
The quotient of these (which we then take the degree $\leq n$ part of)
is as required:
$$\nu(k) \sqcup \Psi(\Lambda_K) \sqcup \left(
1 + \Theta(M,K) - f\left<\frac{1}{2}
\begin{picture}(1,0.7)(-0.25,-0.2)
\qbezier[16](0,-0.1)(0,0.25)(0.25,0.25)
\qbezier[16](0.25,0.25)(0.5,0.25)(0.5,-0.1)
\put(-0.1,-0.4){$k'$}
\put(0.4,-0.4){$k'$}
\end{picture} \ ,\
\psi
\right> + \ldots \right).\eqno{\qed}$$

\section{Changing band self-crossings} 

It seems important to discover a full topological 
formula for $\Theta(M,K)$ (noting
that Lev Rozansky can compute the invariant in all cases \cite{RHP}).
Perhaps 
some enlargement of
Turaev's ``multiplace generalisation of the Seifert matrix'' \cite{Tu}
is the key. Let $\Sigma$ be a Seifert surface for $K$ in $M$ and let
$\{c_1,\ldots,c_{2g}\}$ be a system of curves presenting a basis
for the first homology
of $\Sigma$. Then:

{\bf Problem:} Express $\Theta(M,K)$ in terms of the finite type invariants of degree
$\leq 3$ of the links obtained by pushing the curves $\{c_i\}$ off
$\Sigma$.

We wrap this article up with a step in this direction. 
Let $\Sigma$ be a Seifert surface of genus $2g$
for a knot $(M,K)$. Focus on
a band of this surface (the dashed part, below, follows some arbitrarily
knotty path around the other bands in $M$), 
and consider another knot $(M,K^{{\bowtie}})$ obtained
from $(M,K)$ by a move of the following sort (which we will call
a {\it band twist}): 
\begin{equation}\label{bandtwist}
\epsfxsize 3cm
\epsfbox{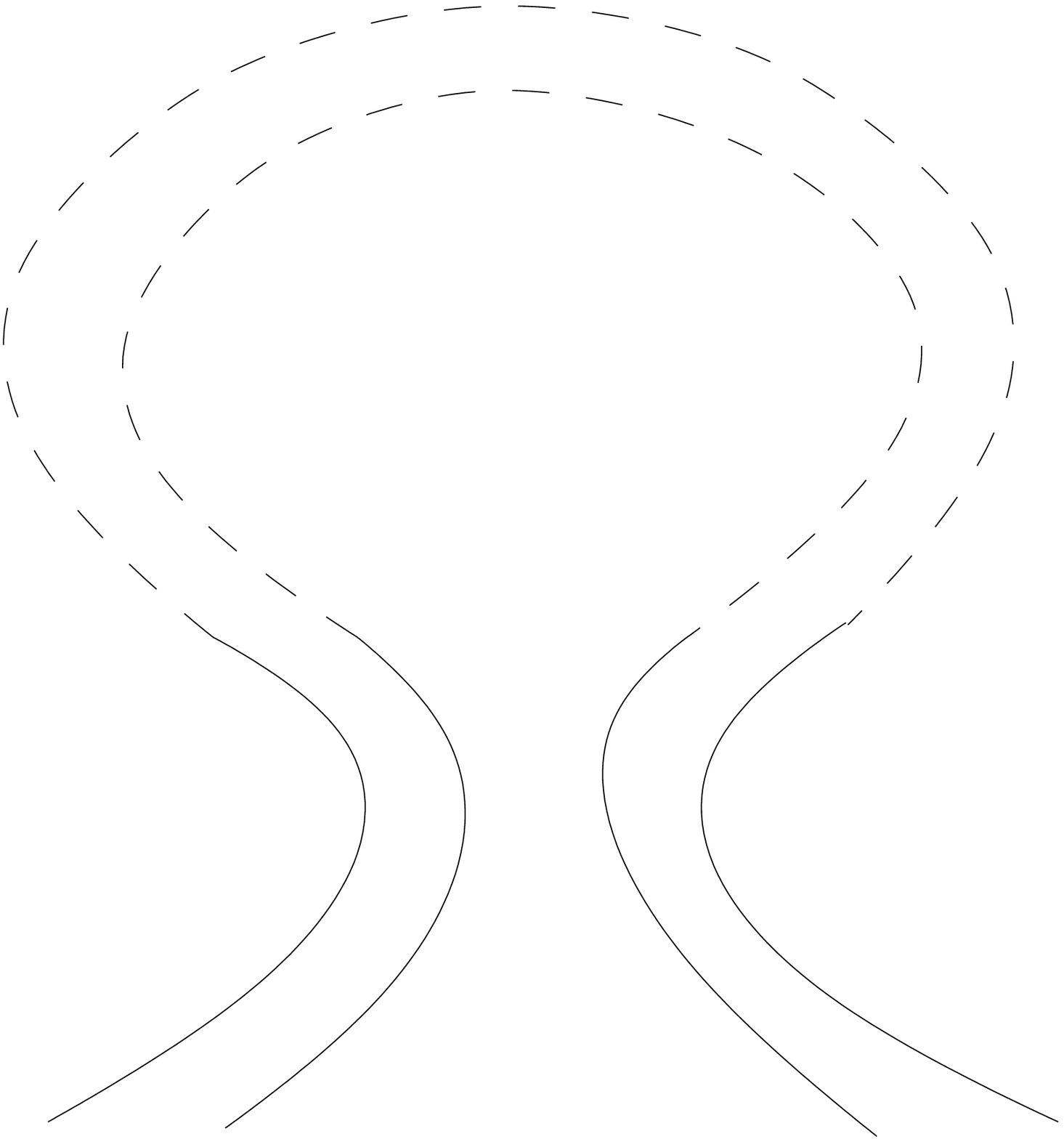}
\hspace{1.15cm}
\begin{picture}(1,1)(0,0)
\put(0,1){$\mapsto$}
\end{picture}
\hspace{0.15cm}
\epsfxsize 3cm
\epsfbox{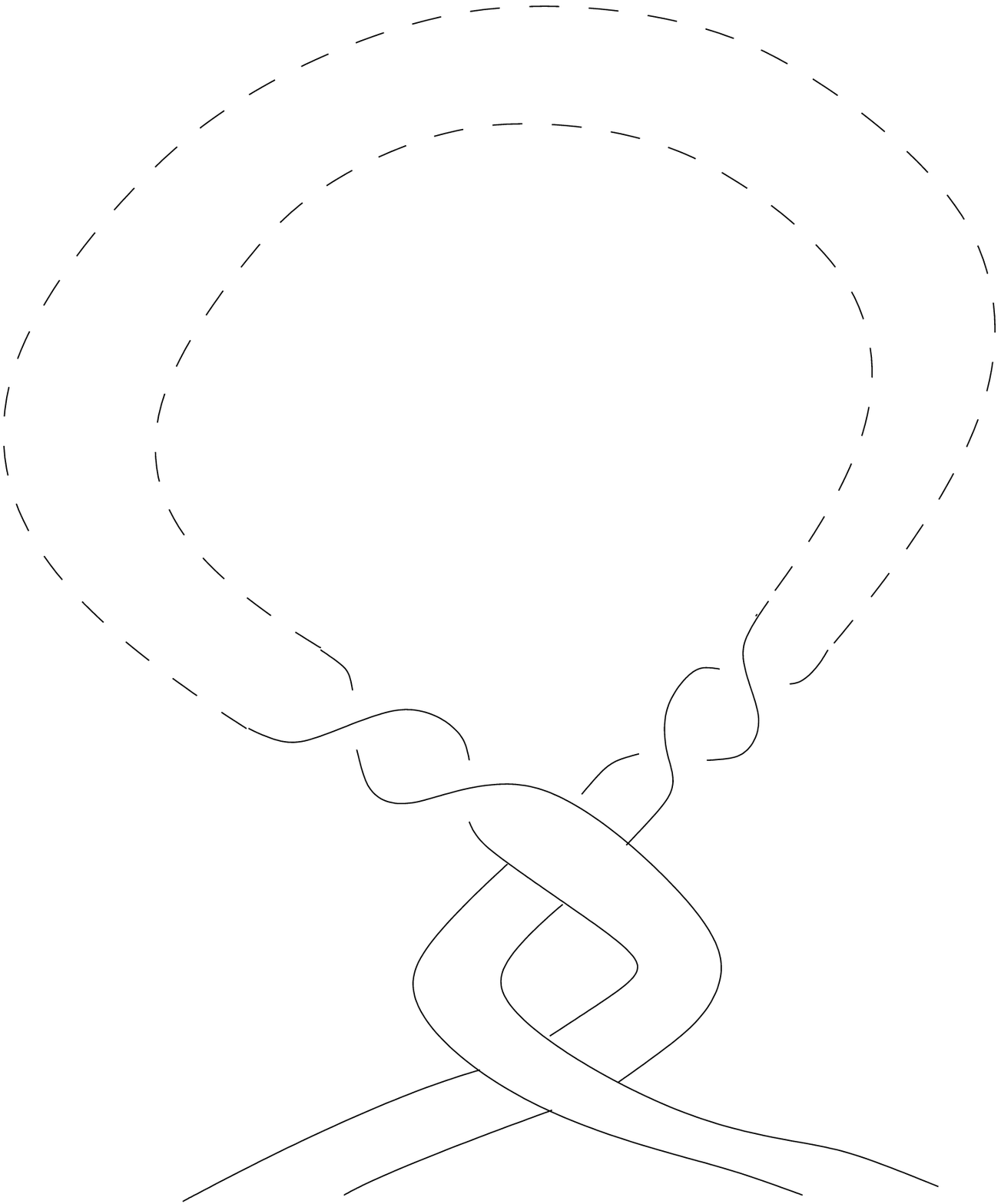}
\end{equation}
We can now apply Theorem \ref{mainth} to 
calculate $\Theta(M,K^{\bowtie}) - \Theta(M,K)$ for the reason
that these two knots are related 
by a surgery on a knot $K'$ which together with
$K$ forms a boundary link in $M$. $K'$ can be taken to be the 
$-1$-framing of the boundary of $\Sigma'$ (see diagram \ref{curves} below).

\begin{equation}\label{curves}
\epsfxsize 3.7cm
\epsfbox{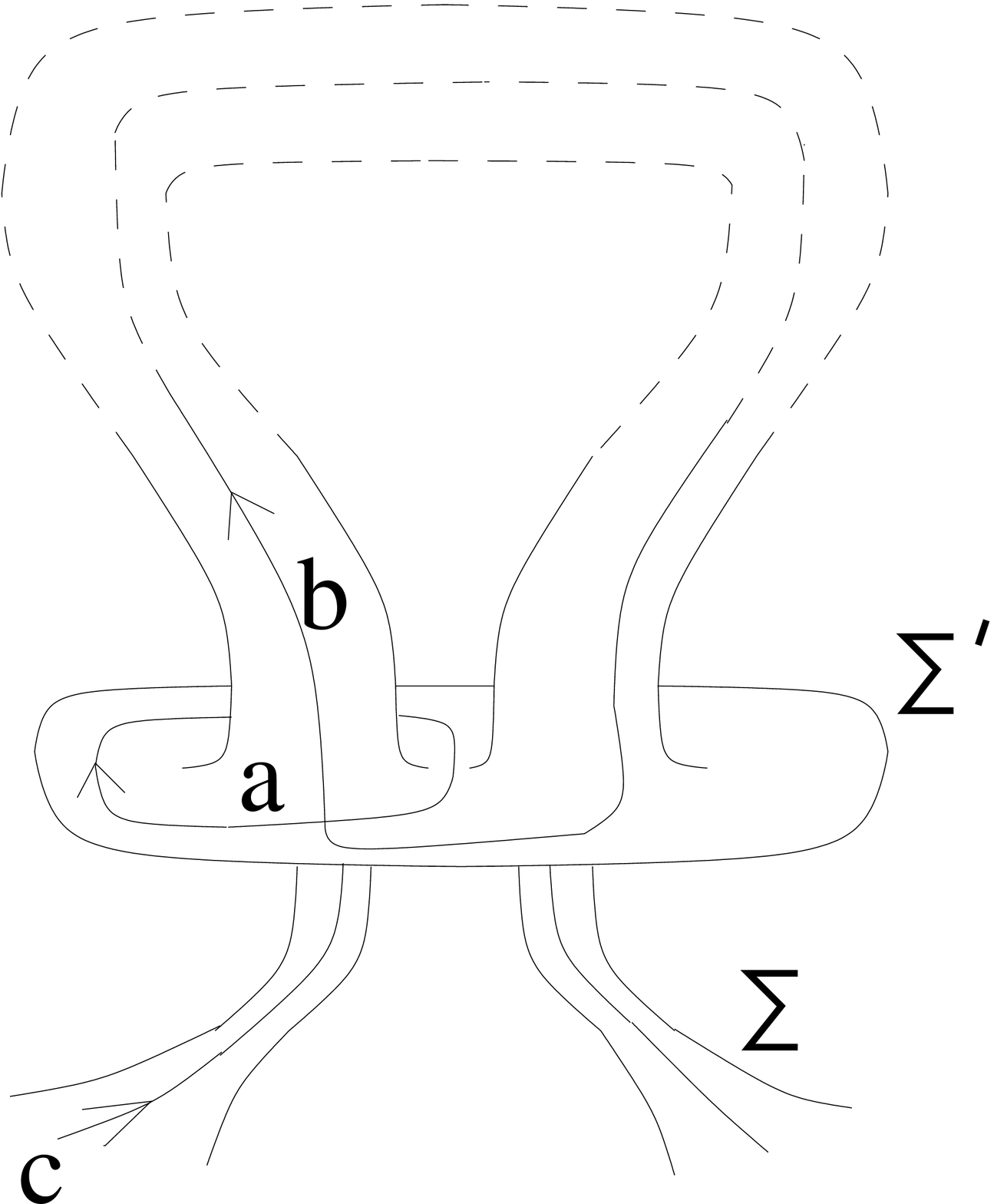}
\end{equation}
First, we choose a convenient
system of curves presenting a basis for the first homology
of $(\Sigma,\Sigma')$.
Let $c_1$ be some curve in $\Sigma$
which traverses the part of the band
shown in the diagram exactly once (shown as $c$).
Complete the choice of the $\{c_i\}$ in any fashion,
so long as the chosen curves lie in the complement in $\Sigma$
of the displayed part of the band.
Let $c'_1$ be $a$ and let $c'_2$ be $b$ (this choice is
fixed by the diagram if 
we demand that lk$(b^+,b)=0$, where $b^+$ is the push-off of $b$). 
The corresponding
$(2g+2)\times(2g+2)$ Seifert matrix can be written:
\[
S = 
\left[
\begin{array}{c|c}
{S_K}  & 
\begin{array}{c}
-1  \\
0  \\
.  \\
.  \\
0 
\end{array}\,\,
{\bf b}^{*} \\
\hline 
\begin{array}{c}
\begin{array}{ccccc}
-1 & 0 & . & . & 0 
\end{array} \\
{\bf b}
\end{array}
& \hspace{0.30cm} 
\begin{array}{cr}
0 & -1 \\
0 & 0
\end{array}
\end{array}
\right],
\]
where the vector ${\bf {b}}$ 
collects the linking numbers of the curve 
$b$ with the system of curves $\{c_i\}$ on $\Sigma$:
\[
{\bf b} = \{ \mbox{lk}(b,c_1),
\mbox{lk}(b,c_2),
\ldots,
\mbox{lk}(b,c_{2g}) \}. \]

Substituting this Seifert matrix into the surgery formula
(Theorem \ref{mainth}), one arrives (via some involved matrix algebra)
at the following statement. Let ${\bf p}$ denote the projector
$\overbrace{(1,0,\ldots,0)}^{2g}$, and let $\left<k\right>$ 
(resp. $\left<k'\right>$) denote 
$e^{k/2}-e^{-k/2}$ (resp. $e^{{k'}/{2}}-e^{-{k'}/{2}}$).
\begin{thm}
Letting $f=-1$ for a positive twist (as
shown earlier), and $f=+1$ for a negative twist:\vspace{-0.2cm}
\[
\Theta(M,K^{\bowtie})
=
\Theta(M,K) - f
\left<
\frac{1}{2}
\begin{picture}(1,0.7)(-0.25,-0.2)
\qbezier[16](0,-0.1)(0,0.25)(0.25,0.25)
\qbezier[16](0.25,0.25)(0.5,0.25)(0.5,-0.1)
\put(-0.1,-0.4){$k'$}
\put(0.4,-0.4){$k'$}
\end{picture} \ ,\
\Psi\left(Z\right)^{-1/2} \right>,
\]
where $Z$ is the matrix
\[
I\ -\ 
\Lambda_{S_K}^{-1}\,
\mbox{$\left<k'\right>$}
\left(e^{k'/2} {\bf b}^* {\bf p} - 
e^{-k'/2} {\bf p}^* {\bf b} \right)
\left<k\right>.
\]
\end{thm}

This is technical, yes, but straightforward to compute. Observe
that it has the expected rationality structure, with the Alexander
polynomial of $K^{\bowtie}$ appearing in the denominators of the labels.

There is
an obvious extension to the case of a twist of $2\pi p$ radians.

\section*{Appendix: Multiplicative corrections}
\setlength{\unitlength}{0.75cm}

Frequently we encounter something like the following problem:
we are presented an element of ${\cal A}(\ostar_{\{k,k'\}})$
as (the projection from ${\cal A}(\star_{\{k,k'\}})$ of)
an exponential (w.r.t. the ``disjoint-union product'')
of some series of connected 
``symmetrised'' diagrams. For example:
\[
\mbox{exp}_\sqcup(W+S_1+S_2+\ldots) \in {\cal A}(\ostar_{\{k,k'\}}),
\] 
where $W$ is a series of wheels all of whose legs are labelled
$k$ and 
$S_i$ is a series of connected 
fragments of generating diagrams, each with exactly 
$i$ trivalent vertices
 (see Definition
\ref{fragment}).
We are then asked to
multiply it by some given factor
using the ``connect-sum multiplication'', 
and then re-express the result
in some convenient ``symmetrised''
form. We will examine, as a guiding example, the following expression
(note that the error term on the LHS below is {\it different} to the
error term on the RHS):

\begin{lem}
\setlength{\unitlength}{0.8cm}
\begin{eqnarray*}
\lefteqn{\mbox{exp}_{\#}
\left(
\frac{f}{2}
\begin{picture}(1.7,0.7)(-0.7,-0.255)
\qbezier[10](0,-0.1)(0,0.25)(0.25,0.25)
\qbezier[10](0.25,0.25)(0.5,0.25)(0.5,-0.1)
\put(-0.5,-0.1){\vector(1,0){1.5}}
\end{picture} 
\right) 
\#_{\{k'\}} \mbox{exp}_\sqcup(W+S_1+S_2+\ldots)} \\
& = &
\mbox{exp}_{\sqcup}\left(
\frac{f}{2}
\begin{picture}(1,0.7)(-0.25,-0.2)
\qbezier[10](0,-0.1)(0,0.25)(0.25,0.25)
\qbezier[10](0.25,0.25)(0.5,0.25)(0.5,-0.1)
\put(-0.2,-0.5){$k'$}
\put(0.4,-0.5){$k'$}
\end{picture}\ +
W + S_1 + S_2 + \ldots \right). 
\end{eqnarray*}
\end{lem}
\begin{proof}
There are brute force ways to do this, but an elegant approach
is available: 
use the {\it wheeling isomorphism}. This theory was conjectured
by Bar-Natan, Garoufalidis, Rozansky and Thurston,
conjectures which have since been proved, including a Kontsevich
integral based theory due to Bar-Natan, Le
and Thurston \cite{BLT}.

All we need to know about
this beautiful game here is the following: there exists an element 
$\nu(k)\in {\cal A}(\star_k)$, which is of the form $\mbox{exp}_{\sqcup}(W)$,
for $W$ a series of connected {\it wheels} (1-loop diagrams), with the
property that if $X\in {\cal A}(\star_{x})$, and if $Y\in 
 {\cal A}(\ostar_{x,y})$ then
\begin{equation}\label{BLT}
\widehat{\nu(x)}(X) \#_{\{x\}} \widehat{\nu(x)}(Y) =
\widehat{\nu(x)}(X\sqcup Y),
\end{equation}
where $\widehat{\nu(x)}(X)$, for example,
is the operation of joining, in all ways,
{\it all}
legs of some diagram appearing in $\nu(x)$ to {\it some} of the 
$x$-labelled legs of some diagram appearing in $X$,
bilinearly extended to all diagrams in those series. 
(At the end of this Appendix we remind why
the BLT theory applies in this generality.)

To apply this here, note that (recalling that $\kappa$ denotes the coefficient 
of the wheel with two legs in $\nu(x)$):
\begin{equation}\label{eqna}
\frac{f}{2}
\ \begin{picture}(1.7,0.7)(-0.7,-0.1)
\qbezier[14](0,-0.1)(0,0.25)(0.25,0.25)
\qbezier[14](0.25,0.25)(0.5,0.25)(0.5,-0.1)
\put(-0.5,-0.1){\vector(1,0){1.5}}
\put(0.65,0.1){$x$}
\end{picture}\ \ 
= 
\widehat{
\nu(x)}
\left( 
\frac{f}{2}
\begin{picture}(1,0.7)(-0.25,-0.2)
\qbezier[12](0,-0.1)(0,0.25)(0.25,0.25)
\qbezier[12](0.25,0.25)(0.5,0.25)(0.5,-0.1)
\put(-0.1,-0.35){$x$}
\put(0.4,-0.35){$x$}
\end{picture} 
-
f
\kappa
\thetdiagsm
\right)
,
\end{equation}
and also that
\begin{equation}\label{eqnb}
\mbox{exp}_\sqcup(W+S_1+S_2+\ldots) =
\widehat{\nu(k')}
\left(
\mbox{exp}_\sqcup(W+S_1+S_2+\ldots)
\right),
\end{equation}
because any way of gluing a wheel to the legs of 
a fragment of a generating diagram will result in
a fragment of a generating diagram with
at least three trivalent vertices.
So, taking the $\#$-exponential of (\ref{eqna}), and then
$\#$-multiplying it into (\ref{eqnb}) gives:
\begin{eqnarray*}
\lefteqn{\widehat{\nu(k')}
\left(
\mbox{exp}_{\sqcup}\left(
\frac{f}{2}
\begin{picture}(1,0.7)(-0.25,-0.2)
\qbezier[12](0,-0.1)(0,0.25)(0.25,0.25)
\qbezier[12](0.25,0.25)(0.5,0.25)(0.5,-0.1)
\put(-0.1,-0.35){$k'$}
\put(0.4,-0.35){$k'$}
\end{picture} 
-
f
\kappa
\thetdiagsm
+
W+S_1 + S_2 + \ldots
\right)\right)}. \\
& = &
\left(
\mbox{exp}_{\sqcup}\left(
\frac{f}{2}
\begin{picture}(1,0.7)(-0.25,-0.2)
\qbezier[12](0,-0.1)(0,0.25)(0.25,0.25)
\qbezier[12](0.25,0.25)(0.5,0.25)(0.5,-0.1)
\put(-0.1,-0.35){$k'$}
\put(0.4,-0.35){$k'$}
\end{picture} 
-
f
\kappa
\thetdiagsm
+
f
\kappa
\thetdiagsm
+
W+S_1 + S_2 + \ldots
\right)\right), 
\end{eqnarray*}
as required.
\end{proof}

We finish with a remark reminding why the BLT theory includes equation \ref{BLT}.
Following \cite{BLT}:
\[
\hat{Z}\left(
\begin{picture}(1,1)(0.1,0.15)
\qbezier(0.5,0.5)(0.75,0.5)(0.75,0.3)
\qbezier(0.6,0.1)(0.75,0.1)(0.75,0.3)
\qbezier(0.5,0.5)(0.25,0.5)(0.25,0.3)
\qbezier(0.4,0.1)(0.25,0.1)(0.25,0.3)
\qbezier(0.5,-0.5)(0.5,0.2)(0.5,0.4)
\qbezier(0.5,0.6)(0.5,0.7)(0.5,1)
\put(0.275,0.4){\vector(0,1){0.01}}
\put(0.5,-0.3){\vector(0,1){0.05}}
\put(0.6,0.7){$x$}
\put(0.8,-0.0){$k$}
\end{picture}
\right) = 
\mbox{exp}_{\sqcup}\left(
\frac{f}{2}
\begin{picture}(1,0.7)(-0.25,-0.2)
\qbezier[12](0,-0.1)(0,0.25)(0.25,0.25)
\qbezier[12](0.25,0.25)(0.5,0.25)(0.5,-0.1)
\put(-0.15,-0.35){$x$}
\put(0.4,-0.35){$k$}
\end{picture}\, \right) \sqcup \nu(k)
\in {\cal A}(\star_{\{x\}}\,\ostar_{\{k\}}).\ 
\]
We may ``double'' the component $k$ in 2 different ways: by taking its parallel
or by vertically composing 2 copies of the long Hopf link (``1+1=2''). The functorial
properties of $\hat{Z}$ lead to a corresponding algebraic equation in 
${\cal A}(\star_{\{x\}}\,\ostar_{\{k_1,k_2\}}).$ Pairing with $X(k_1)\sqcup
Y(k_2,y)$ sends this equation to equation \ref{BLT}. The unusual point that
must be observed is that this pairing is still well-defined (despite the
$y$-labelled legs)
because after pairing a ``link relation'' in $k_2$ may ``sweep'' 
across $Y$ to become a 
link relation in $y$, which we are modding out by.

\Addresses\recd


\begin{thebibliography}

\bibitem{AM}{\bf S Akbulut}, {\bf D McCarthy}, {\it Casson's invariant
for oriented homology 3-spheres, an exposition}, Mathematical notes, 36,
Princeton Univ. Press

\bibitem{BG}{\bf D Bar-Natan}, {\bf S Garoufalidis}, {\it On the 
Melvin-Morton-Rozansky conjecture}, Invent. Math. {125} (1996) 103-133

\bibitem{BGRT}{\bf D Bar-Natan}, {\bf S Garoufalidis}, {\bf L Rozansky},
{\bf D Thurston}, {\it The Aarhus integral of rational homology
3-spheres II: Invariance and universality}, to appear in Selecta Math.
{\tt arXiv:math.QA/9801049}

\bibitem{BGRTW}{\bf D Bar-Natan}, {\bf S Garoufalidis}, {\bf L Rozansky},
{\bf D Thurston}, {\it Wheels, wheeling, and the Kontsevich integral
of the unknot}, Israel J. Math. {119} (2000) 217--237

\bibitem{BLT}{\bf D Bar-Natan}, {\bf T Le}, {\bf D Thurston}, {\it
Two applications of elementary knot theory to Lie algebras and Vassiliev
invariants}, {\tt arXiv:math.GT/0204311}

\bibitem{GGP} {\bf S Garoufalidis}, {\bf M Goussarov}, {\bf M Polyak},
{\it Calculus of clovers and finite type invariants of 3-manifolds},
Geom. Topol. {5} (2001) 75--108

\bibitem{GK} {\bf S Garoufalidis}, {\bf A Kricker}, {\it A rational
noncommutative invariant of boundary links}, 
{\tt arXiv:math.GT/0105028}

\bibitem{GL} {\bf S Garoufalidis}, {\bf J Levine}, {\it Analytic
invariants of boundary links}, to appear in Proceedings of Knots 2000,
J. Knot Th. Ramifications

\bibitem{GR} {\bf S Garoufalidis}, {\bf L Rozansky}, {\it The loop
expansion of the Kontsevich integral, abelian invariants of knots and
S-equivalence}, {\tt arXiv:math.GT/0003187}

\bibitem{K}{\bf A Kricker}, {\it The lines of the Kontsevich
invariant and Rozansky's rationality conjecture}, 
{\tt arXiv:math.GT/0005284}

\bibitem{KS}{\bf A Kricker}, {\it Calculations of the Kontsevich integral
using clasper calculus and surgery formulae}, Proceedings of the
46th Japan Topology Symposium, Sapporo, July 1999, (M Morimoto, editor)

\bibitem{LeD}{\bf T\,T\,Q Le}, {\it On the denominators of the Kontsevich 
integral and the universal perturbative invariant of three-manifolds},
Invent. Math. {135} (1999) 689--722

\bibitem{LeG}{\bf T\,T\,Q Le}, {\it On the LMO invariant}, notes
accompanying lectures given at the ``Summer school on quantum
invariants of knots and three-manifolds'', organised by C Lescop,
Grenoble, July 1999

\bibitem{LMO}{\bf T\,T\,Q Le}, {\bf J Murakami}, {\bf T Ohtsuki},
{\it A universal quantum invariant of three-manifolds}, Topology,
{37} (1998) 539--574

\bibitem{MM}{\bf P Melvin}, {\bf H Morton}, {\it The coloured
Jones function}, Comm. Math. Phys. {169} (1995) 501--520

\bibitem{MO}{\bf J Murakami}, {\bf T Ohtsuki}, {\it Topological
quantum field theory for the universal quantum invariant},
Comm. Math. Phys. {188} (1997) 501--520

\bibitem{BPM} {\bf B Patureau-Mirand}, {\it Non-injectivity of
the ``hair'' map},\nl {\tt arXiv:math.GT/0202065}

\bibitem{RozG}{\bf L Rozansky}, {\it A rational generating structure
for finite-type invariants}, notes accompanying lectures given at
the ``Summer school on quantum invariants of knots and three-manifolds'',
organised by C Lescop, Grenoble, July 1999

\bibitem{RMM}{\bf L Rozansky}, {\it The universal $R$-matrix, Burau
representation, and the Melvin-Morton expansion of the coloured
Jones polynomial}, Adv. Math. {134} (1998) 1--31

\bibitem{RHP}{\bf L Rozansky}, {\it Lev Rozansky's Homepage}, Yale
University,\nl {\tt http://www.math.yale.edu/users/rozansky}

\bibitem{Tu}{\bf V Turaev}, {\it Multiplace generalisations of the
Seifert form of a classical knot}, Math. USSR-Sbornik, {44} (1981)
335--361

\end{thebibliography}
\end{document}